\documentclass[12pt,reqno]{amsart}

\usepackage{amsfonts}
\usepackage{eurosym}
\usepackage{amssymb}
\usepackage{amsthm}
\usepackage{amsmath}
\usepackage{amsaddr}
\usepackage{bm}
\usepackage{cite}
\usepackage{mathrsfs}
\usepackage{xcolor}
\usepackage[OT1]{fontenc}
\usepackage[left=2.0cm, right=2.0cm, top=3cm]{geometry}
\usepackage{hyperref}
\hypersetup{colorlinks=true, linkcolor=blue, citecolor=red}

\RequirePackage{times}
\allowdisplaybreaks
\newtheorem{theorem}{Theorem}[section]

\newtheorem{remark}[theorem]{Remark}

\def\n{\textbf{\textit{n}}}
\def\R3{\mathbb{R}^3}
\def\F2o{\overline{F_2}}

\def\d{{\rm d}}
\def \l {\langle}
\def \r {\rangle}
\def\V{{\mathbf{V}}}
\def\T{\mathbb{T}}

\def\H{\mathbf{H}}
\def\W{\mathbf{W}}
\def\A{\mathbf{A}}

\def\uu{\textbf{\textit{u}}}
\def\vv{\textbf{\textit{v}}}
\def\ww{\textbf{\textit{w}}}
\def\ddt{\frac{\d}{\d t}}

\def\J{\mathbf{J}}

\def \HH{\mathbb{H}}
\def\VV{\mathbb{V}}
\def\WW{\mathbb{W}}

\def\P{\mathbb{P}}
\def\div{\mathrm{div}\,}
\def\L2{L^2(\Omega)}
\def\LT{L^2(\mathbb{T}^2)}

\def \au {\rm}
\def \ti {\it}
\def \jou {\rm}
\def \bk {\it}
\def \no#1#2#3 {{\bf #1} (#3), #2.}
\def \eds#1#2#3 {#1, #2, #3.}
\def \nome#1#2 {{\bf #1}, (#2).}

\begin{document}
\title[Well-Posedness of the AGG model]{Well-posedness of the two-dimensional Abels-Garcke-Gr\"{u}n model \\for two-phase flows with unmatched densities}
\author[A. Giorgini]{Andrea Giorgini}

\address{Department of Mathematics \&
Institute for Scientific Computing and Applied Mathematics\\
 Indiana University\\
Bloomington, IN 47405, USA}
\email{agiorgin@iu.edu}
\date{\today}



\begin{abstract}
We study the Abels-Garcke-Gr\"{u}n (AGG) model for a mixture of two viscous incompressible fluids with different densities. The AGG model consists of a Navier-Stokes-Cahn-Hilliard system characterized by a (non-constant) concentration-dependent density and an additional flux term due to interface diffusion. In this paper we address the well-posedness problem in the two-dimensional case. We first prove the existence of local strong solutions in general bounded domains. In the space periodic setting we show that the strong solutions exist globally in time. In both cases we prove the uniqueness and the continuous dependence on the initial data of the strong solutions.
\end{abstract}

\maketitle


\section{Introduction}
In this article we consider the Abels-Garcke-Gr\"{u}n (AGG) model
\begin{equation}
 \label{AGG}
\begin{cases}
\partial_t (\rho(\phi)\uu) + \div \big( \uu \otimes (\rho(\phi) \uu + \widetilde{\J})\big) - \div (\nu(\phi)D\uu) + \nabla P= - \div(\nabla \phi \otimes \nabla \phi)\\
\div \uu=0\\
\partial_t \phi +\uu\cdot \nabla \phi = \Delta \mu\\
\mu= -\Delta \phi+\Psi'(\phi).
\end{cases}
\end{equation}
The AGG system is studied in $\Omega \times (0,T)$, where $\Omega$ is either a bounded domain in $\mathbb{R}^2$ or the two-dimensional torus $\T^2$. The state variables are the volume averaged velocity $\uu=\uu(x,t)$, the pressure of the mixture $P=P(x,t)$, and the difference of the fluids concentrations $\phi=\phi(x,t)$. The symmetric gradient is $D=\frac12 (\nabla +\nabla^T)$. The flux term $\widetilde{\J}$, the mean density $\rho$ and the mean viscosity $\nu$ of the mixture are given by
\begin{equation}
\label{Jrhonu}
\widetilde{\J}= -\frac{\rho_1-\rho_2}{2}\nabla \mu, \quad \rho(\phi)= \rho_1 \frac{1+\phi}{2}+ \rho_2 \frac{1-\phi}{2}, \quad
\nu(\phi)=\nu_1 \frac{1+\phi}{2}+ \nu_2 \frac{1-\phi}{2},
\end{equation}
where $\rho_1$, $\rho_2$ and $\nu_1$, $\nu_2$ are the homogeneous densities and viscosities of the two fluids.  
The nonlinear function $\Psi$ is the Flory-Huggins potential
\begin{equation}  
\label{Log}
\Psi(s)=\frac{\theta}{2}\bigg[ (1+s)\log(1+s)+(1-s)\log(1-s)\bigg]-\frac{%
\theta_0}{2} s^2, \quad s \in [-1,1],
\end{equation}
where the constant parameters $\theta$ and $\theta_0$ fulfill the conditions $0<\theta<\theta_0$. Notice that \eqref{AGG}$_1$ can be rewritten in the non-conservative form as
\begin{equation}
\label{NS-nc}
\rho(\phi) \partial_t \uu + \rho(\phi) (\uu \cdot \nabla)\uu
-\rho'(\phi) (\nabla \mu\cdot \nabla) \uu - \div (\nu(\phi)D\uu) + \nabla P= - \div(\nabla \phi \otimes \nabla \phi).
\end{equation}
In a bounded domain $\Omega$, the system is subject to the classical boundary conditions
\begin{equation}
\label{AGG-bc}
\uu=\mathbf{0}, \quad \partial_\n\phi=\partial_\n \mu=0 \quad \text{on } \ \partial \Omega \times (0,T),
\end{equation}
where $\n$ is the unit outward normal vector on $\partial \Omega$, and
 $\partial_\n$ denotes the outer normal derivative on $\partial \Omega$.
In the case $\Omega=\T^2$, the state variables satisfy periodic boundary conditions.
In both cases, the system \eqref{AGG} is supplemented with the initial conditions
\begin{equation}
\label{AGG-IC}
\uu(\cdot, 0)=\uu_0, \quad \phi(\cdot, 0)=\phi_0 \quad \text{in } \ \Omega.
\end{equation}
The total energy associated to system \eqref{AGG} is defined as
$$
E(\uu,\phi)= E_{\text{kin}}(\uu, \phi) + E_{\text{free}}(\phi)= 
\int_{\Omega}  \frac12 \rho(\phi) |\uu|^2 \, \d x + \int_{\Omega} \frac12 |\nabla \phi|^2 + \Psi(\phi)   \, \d x,
$$
and the corresponding energy equation reads as
 \begin{equation}
\label{EE}
\ddt E(\uu, \phi) +\int_{\Omega} \nu(\phi) |D \uu|^2 \, \d x+ 
\int_{\Omega} |\nabla \mu|^2 \, \d x =0.
\end{equation}

The Abels-Garcke-Gr\"{u}n system is a fundamental diffuse interface model which describes the motion of two viscous incompressible Newtonian fluids with unmatched densities (i.e. $\rho_1\neq \rho_2$). The model was derived in the seminal paper \cite{AGG2012}. The AGG model is a thermodynamically consistent generalization of the well-known Model H (see \cite{GPV1996} for the derivation and \cite{A2009,GMT2019} for the mathematical analysis). In fact, the classical Navier-Stokes-Cahn-Hilliard system is recovered in the matched density case (i.e. $\rho_1=\rho_2$) since the flux $\widetilde{\J}=\mathbf{0}$ and the density $\rho(\phi)$ is constant. As for the Model H, the fluid mixture in the AGG system is driven by the capillarity forces $-\div(\nabla \phi \otimes \nabla \phi)$ due to the surface tension effect. In addition, a partial diffusive mixing is assumed in the interfacial region, which is modeled by $\Delta \mu$, being the chemical potential $\mu= \frac{\delta E_{\text{free}}(\phi)}{\delta \phi}$. The specificity of the AGG model lies in the presence of the flux term $\widetilde{\J}$. In contrast to the one-phase flow, the (average) density $\rho(\phi)$ in \eqref{AGG} does not satisfy the continuity equation with respect to the flux associated with the velocity $\uu$. Instead, the density $\rho(\phi)$ satisfies the continuity equation with a flux given by the sum of the transport term $\rho(\phi)\uu$ and the term $\widetilde{\J}$. The latter is due to the diffusion of the concentration in the unmatched densities case.
For the connection of the AGG model with the classical sharp interface two-phase problem and the so-called sharp interface limit, 
we refer the reader to the review in \cite{AG2018}. 
It is important to mention that the theory of diffuse interface models for mixtures of fluids has been widely developed in the past decades. Several systems have been proposed to model binary mixtures with non-constant density in view of their applications in engineering and physics. We mention the models derived in \cite{B2002, DSS2007, GKL2018, HMR2012, LT1998, SSBZ2017} and the theoretical analysis achieved in \cite{A2009-2, A2012, B2001, GT2020, KZ2015}.

The mathematical analysis of the AGG system has been focused so far on the existence of weak solutions in two and three dimensional bounded domains. More precisely, global solutions with finite energy for the system \eqref{AGG} with boundary and initial conditions \eqref{AGG-bc}-\eqref{AGG-IC} were proven in \cite{ADG2013} and \cite{ADG2013-2}.
In the former the mobility coefficient $m(\phi)$ is non-constant and strictly positive, whereas in the latter $m(\phi)$ is degenerate\footnote{In comparison with  \cite{ADG2013,ADG2013-2}, we have set the mobility $m(\phi)$ and the energy coefficient $a(\phi)$ equal to one  in \eqref{AGG}.}. More recently, the existence of global weak solutions have been generalized in \cite{AB2018} for viscous non-Newtonian binary fluids and in \cite{GalGW2019} for the case of dynamic boundary conditions describing moving contact lines. Furthermore, non-local variants of the AGG system have been investigated in \cite{AT2020} and in \cite{Frigeri2016}, where the gradient term $\frac12 | \nabla \phi|^2$ in the  local free energy $ E_{\text{free}}(\phi)$ has been replaced by different non-local operators. 

The aim of this contribution is to present the first well-posedness result for the AGG model. In our analysis we show existence, uniqueness and continuous dependence on the initial data of the strong solutions in the two-dimensional case. In comparison with the notion of weak solutions studied in all the previous works on the AGG model, such strong solutions are more regular and solve the system \eqref{AGG} pointwise almost everywhere. These solutions depart from initial data $\uu_0 \in \H^1(\Omega)$ with  $\div \uu_0=0$, and $\phi_0 \in H^2(\Omega)$ such that $-1\leq \phi_0(x)\leq 1$ in $\Omega$ and $-\Delta \phi_0+\Psi'(\phi_0)\in H^1(\Omega)$, which satisfy suitable boundary or periodic conditions.
We first prove the existence of local-in-time strong solutions in a general bounded domain (see Theorem \ref{mr1}).  
Our proof relies on the existence of suitable (global) approximate solutions to system \eqref{AGG} constructed through a semi-Galerkin formulation. In this framework the modified Navier-Stokes equations \eqref{AGG}$_1$-\eqref{AGG}$_2$ are solved in finite-dimensional (spacial) spaces, whereas the convective Cahn-Hilliard system \eqref{AGG}$_3$-\eqref{AGG}$_4$ is fully solved (i.e. not approximated). The advantage of this approach is that the approximate velocity fields $\uu_m$ is regular in the space variable, and the approximate concentrations $\phi_m$ take values in the physical interval $[-1,1]$ which, in turn, ensures that 
$\rho'(\phi)=\frac{\rho_1-\rho_2}{2}$\footnote{A different approximation leading to a concentration $\phi_m$ with values outside the interval $[-1,1]$ may need a suitable extension of $\rho(\cdot)$ outside the interval $[-1,1]$, and, in general, it may happen that $\rho'(\phi)\neq \frac{\rho_1-\rho_2}{2}$.}. 
It is worth pointing out that our strategy entirely exploit the regularity properties of the Cahn-Hilliard equation with logarithmic potential in two dimensions. More precisely, the control of $\Psi''(\phi)$ in $L^p$ spaces (available in the two dimensional setting) allows us to recover the time continuity of the chemical potential $\mu$, which is needed to solve the approximated problem. 
Once the existence of the approximate solutions is shown, we employ the energy method to deduce uniform estimates and the necessary compactness to obtain the existence of a local solution to \eqref{AGG}.
Next, in the periodic boundary setting we demonstrate that the strong solutions exist globally in time (see Theorem \ref{mr2}). 
The key observation to obtain the propagation of regularity for all times is that global-in-time higher-order estimates for the full system as in  \cite{GMT2019,GT2020} are out of reach due to the presence of the nonlinear term $(\nabla \mu \cdot \nabla) \uu$ (cf. the term $I_3$ in \eqref{I3}). Notice that, since $\nabla \mu$ belongs to 
$L^2(0,T;L^2(\T^2))$ (cf. \eqref{EE}), $\nabla \mu$ has a lower regularity than $\uu$. Therefore, the idea is to split the argument by first improving the regularity of the concentration $\phi$ relying on the energy estimates obtained from \eqref{EE}, and then showing more regularity properties for the velocity field. A similar idea was used in \cite{A2009} for the Model H. However, the argument in \cite[Lemma 3]{A2009} is based on the integrability properties of $\partial_t \uu$ or the fractional in time regularity of $\uu$, which are not known for the weak solutions to \eqref{AGG}. Nevertheless, it is possible to overcome this issue by exploiting the fine structure of the incompressible Navier-Stokes equations in the periodic setting. The crucial term involving the time derivative of the velocity is rewritten in \eqref{Key-term} in such a way that the highest space derivative acting on the velocity is of order one, and boundary terms do not appear when integrating by parts. Such technique requires
an estimate of the pressure $P$ in $L^2$, which is deduced from the incompressibility condition \eqref{AGG}$_2$ and the crucial estimate 
\eqref{phiH2eT} for the Cahn-Hilliard equation.
Lastly, in both cases (bounded domains and the periodic setting) we show the uniqueness of the strong solutions and their continuous dependence on the initial data.
\medskip

\noindent
\textbf{Plan of the paper.} We report in Section \ref{2} the function spaces  and the notation used in this paper. In Section \ref{MR} we state the main results. Section \ref{Loc-Ex} is devoted to the local existence of strong solutions in bounded domains. In Section \ref{Glo-Ex} we prove the global existence of strong solutions in the space periodic setting. In Section \ref{Un} we address the uniqueness and the continuous dependence on the initial data of the strong solutions.

\section{Preliminaries}
\label{2}
\setcounter{equation}{0}

For a real Banach space $X$, its norm is denoted by $\|\cdot\|_{X}$.
The symbol $\langle\cdot, \cdot\rangle_{X',X}$ stands for the duality pairing between $X$ and its dual space $X'$. The boldface letter $\bm{X}$ denotes the vectorial space endowed with the product structure.
We assume that $\Omega$ is a bounded domain in $\mathbb{R}^2$ with boundary $\partial \Omega$ of class $C^3$ or the flat torus $\T^2=(\mathbb{R}/ 2\pi \mathbb{Z})^2$. 
We denote the Lebesgue spaces by $L^p(\Omega)$ $(p\geq 1)$  with norms $\|\cdot\|_{L^p(\Omega)}$. The inner product in the Hilbert space $L^2(\Omega)$  is denoted by
$(\cdot, \cdot)$.
For $s \in \mathbb{N}$, $p\geq 1$, $W^{s,p}(\Omega)$
is the Sobolev space with norm $\|\cdot\|_{W^{s,p}(\Omega)}$. If $p=2$, we use the notation $W^{s,p}(\Omega)=H^s(\Omega)$. 
For $s=1$ we denote the duality between $H^1(\Omega)$ and the dual space $(H^1(\Omega))'$ by $\l \cdot, \cdot \r$.
In the case $\Omega=\T^2$, we recall that the functions are characterized by their Fourier expansion
$$
f= \sum_{ k\in \mathbb{Z}^2} \widehat{f}_k \mathrm{e}^{i k\cdot x}, 
\quad \text{where} \quad \widehat{f}_{-k}=\overline{\widehat{f}_{k}}^c, \quad \widehat{f_k}= \frac{1}{(2\pi)^2} \int_{\T^2} f(x) \mathrm{e}^{- i k\cdot x} \, \d x,
$$
where $\overline{z}^c$ is the complex conjugate of $z \in \mathbb{C}$.
We report that $\big( \sum_{k \in \mathbb{Z}^2} (1+|k|^{2s}) |\widehat{f}_k|^2 \big)^\frac12$ is a norm on $H^s(\T^2)$, $s \in \mathbb{N}$, which is equivalent to the standard norm.
For every $f\in (H^1(\Omega))'$, we denote by $\overline{f}$ the generalized mean value over $\Omega$ defined by
$\overline{f}=|\Omega|^{-1}\l f,1\r$. If $f\in L^1(\Omega)$, then $\overline{f}=|\Omega|^{-1}\int_\Omega f \, \d x$.
By the generalized Poincar\'{e} inequality, there exists a positive constant $C$ such that 
\begin{equation}
\label{normH1-2}
\| f\|_{H^1(\Omega)}\leq C \big(\| \nabla f\|_{L^2(\Omega)}^2+ |\overline{f}|^2\big)^\frac12, \quad \forall \, f \in H^1(\Omega).
\end{equation}
We recall the Ladyzhenskaya, Agmon and Gagliardo-Nirenberg interpolation inequalities in two dimensions
\begin{align}
\label{LADY}
&\| f\|_{L^4(\Omega)}\leq C \|f\|_{L^2(\Omega)}^{\frac12}\|f\|_{H^1(\Omega)}^{\frac12},  &&\forall \, f \in H^1(\Omega),\\
\label{Agmon2d}
&\| f\|_{L^\infty(\Omega)}\leq C \|f\|_{L^2(\Omega)}^{\frac12}\|f\|_{H^2(\Omega)}^{\frac12},  && \forall \, f \in H^2(\Omega),\\
\label{GN-L4}
&\| \nabla f\|_{L^4(\Omega)}\leq C\| f \|_{L^\infty(\Omega)}^\frac12 \| f \|_{H^2(\Omega)}^\frac12, && \forall \, f \in H^2(\Omega),\\
\label{GN-Linf}
&\| \nabla f\|_{L^\infty(\Omega)}\leq C_s  \| f \|_{L^\infty(\Omega)}^\frac{s-2}{2(s-1)}\| f \|_{W^{2,s}(\Omega)}^\frac{s}{2(s-1)}, && \forall \, f \in W^{2,s}(\Omega), \ s>2.
\end{align}
Next, we introduce the Hilbert spaces of solenoidal vector-valued functions.
In the case of a bounded domain $\Omega \subset \mathbb{R}^2$, we define
\begin{align*}
&\H_\sigma=\{ \uu\in \mathbf{L}^2(\Omega): \mathrm{div}\, \uu=0 \ \text{in } \Omega,\ \uu\cdot \n =0\ \text{on}\ \partial \Omega\},\\
& \V_\sigma =\{ \uu\in \mathbf{H}^1(\Omega): \mathrm{div}\, \uu=0 \ \text{in } \Omega,\ \uu=\mathbf{0}\  \text{on}\ \partial \Omega\}.
\end{align*}
We also use $( \cdot ,\cdot )$ and 
$\| \cdot \|_{L^2(\Omega)}$ for
the inner product and the norm in $\H_\sigma$. The space $\V_\sigma$ is endowed with the inner product and norm
$( \uu,\vv )_{\V_\sigma}=
( \nabla \uu,\nabla \vv )$ and  $\|\uu\|_{\V_\sigma}=\| \nabla \uu\|_{L^2(\Omega)}$, respectively.
We report the Korn inequality
\begin{equation}
\label{KORN}
\|\nabla\uu\|_{L^2(\Omega)} \leq \sqrt2\|D\uu\|_{L^2(\Omega)}, \quad \forall \, \uu \in \V_\sigma,
\end{equation}
which implies that $\| D \uu\|_{\L2}$ is a norm on $\V_\sigma$ equivalent to $\| \uu\|_{\V_\sigma}$.
We introduce the space
$\W_\sigma= \mathbf{H}^2(\Omega)\cap \V_\sigma$
with inner product and norm
$ ( \uu,\vv)_{\W_\sigma}=( \A\uu, \A \vv )$ and $\| \uu\|_{\W_\sigma}=\|\A \uu \|$, where $\A=\mathbb{P}(-\Delta)$ is the Stokes operator and $\mathbb{P}$ is the Leray projection from $\mathbf{L}^2(\Omega)$ onto $\H_\sigma$.
We recall that there exists a positive constant $C>0$ such that
\begin{equation}
\label{H2equiv}
 \| \uu\|_{H^2(\Omega)}\leq C\| \uu\|_{\W_\sigma}, \quad \forall \, \uu\in \W_\sigma.
\end{equation}
In the space periodic case $\Omega=\T^2$, 
we define\footnote{In contrast to the classical periodic setting for the incompressible Navier-Stokes (cf. \cite{Temam}), we do not require that $\widehat{\uu}_0=0$ in the definition of $\HH_\sigma$ and $\VV_\sigma$. This is due to the fact that $\overline{\uu}=\frac{1}{|\Omega|}\int_{\Omega} \uu \, \d x$ is not conserved by the flow of \eqref{AGG}.}
\begin{align*}
&\HH_\sigma=\{ \uu\in \mathbf{L}^2(\T^2): \widehat{\uu}_k \cdot k=0 \quad \forall \, k \in \mathbb{Z}^2 \}, \quad  \VV_\sigma=\H^1(\T^2)\cap \HH_\sigma,
\quad \WW_\sigma= \H^2(\T^2)\cap \HH_\sigma,
\end{align*}
which are endowed with the norms $\| \uu\|_{\HH_\sigma}= \| \uu\|_{L^2(\T^2)}$, $\| \uu\|_{\VV_\sigma}= \| \uu\|_{H^1(\T^2)}$, and 
$\| \uu\|_{\WW_\sigma}= \| \uu\|_{H^2(\T^2)}$.
Since 
\begin{equation}
\label{KORN2} 
\|\nabla\uu\|_{L^2(\T^2)} \leq \sqrt2\|D\uu\|_{L^2(\T^2)}, \quad \forall \, \uu \in \VV_\sigma,
\end{equation}
it follows that $(\| \uu\|_{\LT}^2+ \| D \uu\|_{\LT}^2 )^\frac12$ is a norm on $\VV_\sigma$, which is equivalent to $\| \uu\|_{\VV_\sigma}$.
We recall that
\begin{equation}
\label{H2equiv-T}
 \| \uu\|_{H^2(\T^2)}\leq C(\| \uu\|_{\LT}+\| \Delta \uu\|_{\L2}), \quad \forall \, \uu\in \WW_\sigma.
\end{equation}

Throughout this paper we make use of the following notation:
\\

\noindent
$\bullet$ We define the positive constants
$$
\rho_\ast=\min \lbrace \rho_1,\rho_2\rbrace,
\quad \rho^\ast=\max \lbrace \rho_1,\rho_2\rbrace,
\quad \nu_\ast =\min \lbrace \nu_1,\nu_2 \rbrace,
\quad \nu^\ast =\max \lbrace \nu_1,\nu_2 \rbrace.
$$

\noindent
$\bullet$ We denote the convex part of the Flory-Higgins potential by $F$, namely
$$
F(s)=\frac{\theta}{2}\big[ (1+s)\log(1+s)+(1-s)\log(1-s)\big], \quad s \in [-1,1].
$$

\noindent
$\bullet$ The symbol $C$ denotes a generic positive constant whose value may change from line to line. The specific value depends on the domain $\Omega$ and the parameters of the system, such as $\rho_\ast$, $\rho^\ast$, $\nu_\ast$, $\nu^\ast$, $\theta$ and $\theta_0$. Further dependencies will be  specified when necessary.

\section{Main results}
\label{MR}
\setcounter{equation}{0}

In this section we formulate the main results of this paper.
We state the local well-posedness of the strong solutions to system \eqref{AGG} in a bounded domain $\Omega \subset \mathbb{R}^2$ subject to the boundary conditions \eqref{AGG-bc}.

\begin{theorem}
\label{mr1}
Let $\Omega$ be a bounded domain of class $C^3$ in $\mathbb{R}^2$. Assume that $\uu_0 \in \V_\sigma$ and $\phi_0 \in H^2(\Omega)$ such that $\| \phi_0\|_{L^\infty(\Omega)}\leq 1$, $|\overline{\phi_0}|<1$, $\mu_0= -\Delta \phi_0+ \Psi'(\phi_0) \in H^1(\Omega)$, and $\partial_\n \phi_0=0$ on $\partial \Omega$. Then, there exist $T_0>0$, depending on the norms of the initial data, and a unique strong solution $(\uu, P, \phi)$ to system \eqref{AGG} subject to \eqref{AGG-bc}-\eqref{AGG-IC} on $(0,T_0)$ in the following sense:
\begin{itemize}
\item[(i)] The solution $(\uu, P, \phi)$ satisfies the properties
\begin{align}
\label{reg-SS}
\begin{split}
&\uu \in C([0,T_0]; \V_\sigma) \cap L^2(0,T_0;\W_\sigma)\cap W^{1,2}(0,T_0;\H_\sigma),\\
&P \in L^2(0,T_0;H^1(\Omega)),\\
&\phi \in L^\infty(0,T_0;H^3(\Omega)), \
\partial_t \phi \in L^\infty(0,T_0;(H^1(\Omega))')\cap L^2(0,T_0;H^1(\Omega)),\\
&\phi \in L^\infty(\Omega\times (0,T_0)) : |\phi(x,t)|<1 \ \text{a.e. in } \  \Omega\times(0,T_0),\\
&\mu \in C([0,T_0];H^1(\Omega))\cap L^2(0,T_0;H^3(\Omega))\cap W^{1,2}(0,T_0;(H^1(\Omega))'), \\
&F'(\phi), F''(\phi), F'''(\phi) \in L^\infty(0,T_0;L^p(\Omega)), \ \forall \, p \in [2,\infty).
\end{split}
\end{align}

\item[(ii)] The solution $(\uu, P, \phi)$ fulfills the system \eqref{AGG} almost everywhere in $\Omega \times (0,T_0)$ and the boundary conditions $\partial_\n \phi=\partial_\n \mu=0$ almost everywhere in $\partial \Omega \times (0,T_0)$. 

\item[(iii)] The solution $(\uu, P, \phi)$ is such that $\uu(\cdot, 0)=\uu_0$ and $\phi(\cdot, 0)=\phi_0$ in $\Omega$. Moreover, $(\uu,\phi)$ depends continuously on the initial data in $\H_\sigma\times H^1(\Omega)$ on $[0,T_0]$.
\end{itemize}
\end{theorem}

In the space periodic setting we establish the global well-posedness of the strong solutions. 

\begin{theorem}
\label{mr2}
Let $\Omega= \mathbb{T}^2$. 
Assume that $\uu_0 \in \VV_{\sigma}$ and $\phi_0 \in H^2(\T^2)$ such that $\| \phi_0\|_{L^\infty(\T^2)}\leq 1$, $|\overline{\phi_0}|<1$, 
$\mu_0=-\Delta \phi_0+ \Psi'(\phi_0) \in H^1(\T^2)$. Then, there exists a unique global strong solution $(\uu, P, \phi)$ to system \eqref{AGG} with periodic boundary conditions and initial conditions \eqref{AGG-IC} in the following sense:
\begin{itemize}
\item[(i)] For all $T>0$, the solution $(\uu, P, \phi)$ is such that
\begin{align}
\label{reg-SSP}
\begin{split}
&\uu \in C([0,T]; \VV_\sigma) \cap L^2(0,T;\WW_\sigma)\cap W^{1,2}(0,T;\HH_\sigma),\\
&P \in L^2(0,T;H^1(\T^2)),\\
&\phi \in L^\infty(0,T;H^3(\T^2)), \
\partial_t \phi \in L^\infty(0,T;(H^1(\T^2))')\cap L^2(0,T;H^1(\T^2)),\\
&\phi \in L^\infty(\T^2\times (0,T)) : |\phi(x,t)|<1 \ \text{a.e. in } \  \T^2\times(0,T),\\
&\mu \in C([0,T];H^1(\T^2))\cap L^2(0,T;H^3(\T^2))\cap W^{1,2}(0,T;(H^1(\T^2))'), \\
&F'(\phi), F''(\phi), F'''(\phi) \in L^\infty(0,T;L^p(\T^2)),\ \forall \, p \in [2,\infty).
\end{split}
\end{align}

\item[(ii)] The solution $(\uu, P, \phi)$ satisfies the system \eqref{AGG} almost everywhere in $\T^2 \times (0,T)$.

\item[(iii)] The solution $(\uu, P, \phi)$ fulfills $\uu(\cdot, 0)=\uu_0$ and $\phi(\cdot, 0)=\phi_0$ in $\T^2$. In addition,  for all $T>0$, $(\uu,\phi)$ depends continuously on the initial data in $\HH_\sigma \times H^1(\T^2)$ on $[0,T]$.
\end{itemize}
\end{theorem}


\section{Proof of Theorem \ref{mr1}: Local Existence in Bounded Domains}
\label{Loc-Ex}
\setcounter{equation}{0}

In this section, we prove the existence of local strong solutions to system \eqref{AGG} with boundary and initial conditions \eqref{AGG-bc}-\eqref{AGG-IC} in a bounded domain $\Omega$ in $\mathbb{R}^2$. We first present the semi-Galerkin approximation scheme, then prove the solvability of the approximated system through a fixed point argument, and finally carry out the uniform estimates of the approximate solutions which allow the passage to the limit in the approximate formulation. 

\subsection{Definition of the Approximate Problem}
We consider the family of eigenfunctions  $\lbrace \ww_j\rbrace_{j=1}^\infty$ and eigenvalues $\lbrace \lambda_j\rbrace_{j=1}^\infty$ of the Stokes operator $\A$. For any integer $m\geq 1$, we define the finite-dimensional subspaces of $\V_\sigma$ by $\V_m= \text{span}\lbrace \ww_1,...,\ww_m\rbrace$.
We denote by $\P_m$ the orthogonal projection on $\V_m$ with respect to the inner product in $\H_\sigma$. Since $\Omega$ is of class $C^3$, it follows that $\ww_j \in \mathbf{H}^3(\Omega)\cap \V_\sigma$ for all $j\in \mathbb{N}$. Moreover, we report the inverse Sobolev embedding inequalities in $\V_m$
\begin{equation}
\label{Rev-SI}
\| \vv\|_{H^1(\Omega)}\leq C_m \| \vv\|_{\L2},\quad
\| \vv\|_{H^2(\Omega)}\leq C_m \| \vv\|_{\L2},\quad
\| \vv \|_{H^3(\Omega)}\leq C_m \| \vv\|_{\L2}, \quad 
\forall \, \vv \in \V_m.
\end{equation}

Let us fix $T>0$. For any $m \in \mathbb{N}$, we determine the approximate solution $(\uu_m, \phi_m)$ to the system \eqref{AGG} with boundary and initial conditions \eqref{AGG-bc}-\eqref{AGG-IC} as follows:
\begin{align}
\label{reg-AS}
\begin{split}
&\uu_m \in C^1([0,T];\V_m),\\
&\phi_m \in L^\infty(0,T;W^{2,p}(\Omega)), \
\partial_t \phi_m \in L^\infty(0,T;(H^1(\Omega))')\cap L^2(0,T;H^1(\Omega)),\\
&\phi_m \in L^\infty(\Omega\times (0,T)) : |\phi_m(x,t)|<1 \ \text{a.e. in } \  \Omega\times(0,T),\\
&\mu_m \in C(0,T;H^1(\Omega))\cap L^2(0,T;H^3(\Omega))\cap W^{1,2}(0,T;(H^1(\Omega))'), \\
&F''(\phi_m)\in L^\infty(0,T;L^p(\Omega)),
\end{split}
\end{align}
for all $p \in [2,\infty)$, such that
\begin{align}
\label{NS-w}
\begin{split}
&\begin{aligned}[t]
(\rho(\phi_m) \partial_t \uu_m, \ww)&+(\rho(\phi_m)(\uu_m\cdot \nabla)\uu_m,\ww)+(\nu(\phi_m) D \uu_m, \nabla \ww)\\
&-\frac{\rho_1-\rho_2}{2} ( (\nabla \mu_m \cdot \nabla) \uu_m, \ww)= (\nabla \phi_m \otimes \nabla \phi_m, \nabla \ww), 
\end{aligned}
\end{split}
\end{align}
for all $\ww \in \V_m$ and $t \in [0,T]$, and  
\begin{equation}
\label{CH-w}
\begin{cases}
\partial_t \phi_m +\uu_m \cdot \nabla \phi_m = \Delta \mu_m\\
\mu_m= -\Delta \phi_m+\Psi'(\phi_m)
\end{cases}
\quad \text{a.e. in } \ \Omega \times (0,T).  
\end{equation}
The approximate solution $(\uu_m,\phi_m)$ satisfies the boundary and initial conditions 
\begin{equation}
\label{bc-AS}
\begin{cases}
\uu_m=\mathbf{0}, \quad \partial_\n \phi_m=\partial_\n \mu_m=0 \quad &\text{on } \partial \Omega \times (0,T),\\
\uu_m(\cdot,0)=\P_m \uu_{0}, \ \phi(\cdot,0)=\phi_{0} \quad &\text{in } \Omega.
\end{cases}
\end{equation}

\subsection{Existence of Approximate Solutions}
We perform a fixed point argument to show the existence of the approximate solutions satisfying \eqref{reg-AS}-\eqref{bc-AS}.
To this aim, we take $\vv \in W^{1,2}(0,T;\V_m)$.
We consider the convective Cahn-Hilliard system
\begin{equation}
\label{cCH}
\begin{cases}
\partial_t \phi_m +\vv \cdot \nabla \phi_m = \Delta \mu_m\\
\mu_m= -\Delta \phi_m+F'(\phi_m)-\theta_0 \phi_m
\end{cases}
\quad \text{in }\ \Omega \times (0,T),
\end{equation}
which is equipped with the boundary and initial conditions
\begin{equation}
\label{cCH-c}
\partial_\n \phi_m=\partial_\n \mu_m=0 \quad \text{on} \ \partial \Omega \times (0,T), \quad \phi_m(\cdot, 0)=\phi_0 \quad \text{in }\ \Omega.
\end{equation}
It is proven in \cite[Theorem 6 and Lemma 3]{A2009} that there exists a unique solution to \eqref{cCH}-\eqref{cCH-c} such that
\begin{equation}
\label{phi-As}
\begin{split}
&\phi_m \in L^\infty(0,T;W^{2,p}(\Omega)), \quad 
\partial_t \phi_m \in L^\infty(0,T;(H^1(\Omega))')\cap L^2(0,T;H^1(\Omega)),\\
&\phi_m \in L^\infty(\Omega\times (0,T)) : |\phi_m(x,t)|<1 \ \text{a.e. in } \  \Omega\times(0,T),\\
&\mu_m \in L^\infty(0,T;H^1(\Omega)),
\end{split}
\end{equation} 
for any $p\in [2,\infty)$. Thanks to \cite[Lemma A.6]{CG2020}, it follows that $F''(\phi)\in L^\infty(0,T;L^p(\Omega))$ for any $p \in [2, \infty)$. In addition, by comparison in \eqref{cCH}$_1$ and \eqref{cCH}$_2$, we infer that $\mu \in L^2(0,T;H^3(\Omega))$ and 
$\partial_t \mu_m\in L^2(0,T;(H^1(\Omega))')$ (see, e.g., \cite[Proof of Theorem 5.1]{G2020}).
Therefore, we have 
\begin{equation}
\label{mu-AS}
\mu_m \in C(0,T;H^1(\Omega))\cap L^2(0,T;H^3(\Omega))\cap W^{1,2}(0,T;(H^1(\Omega))').
\end{equation}
We report the following estimates for the system \eqref{cCH}-\eqref{cCH-c} (see \cite{A2009}, cf. also \cite{CG2020,GGW2018}):
\begin{itemize}
\item[1.] $L^2$ estimate:
\begin{equation}
\label{L2-AS}
\sup_{t \in [0,T]}\| \phi_m(t)\|_{L^2(\Omega)}^2 + \int_0^T \| \Delta \phi_m(\tau)\|_{L^2(\Omega)}^2 \, \d \tau \leq \| \phi_0\|_{L^2(\Omega)}^2 + \frac{\theta_0^2}{2}T.
\end{equation}
\item[2.] Energy estimate:
\begin{equation}
\label{EE-AS}
\begin{split}
\sup_{t\in [0,T]} \int_{\Omega} \frac12 | \nabla\phi_m(t)|^2 &+ F(\phi_m(t)) \, \d x + \frac12 \int_0^T \| \nabla \mu(\tau)\|_{L^2(\Omega)}^2 \, \d \tau \\
&\leq  E_{\text{free}}(\phi_0)+ \frac12 \int_0^T \| \vv(\tau)\|_{L^2(\Omega)}^2 \, \d \tau
+ \frac{\theta_0}{2} \| \phi_0\|_{L^2(\Omega)}^2 + \frac{\theta_0^3}{4}T.
\end{split} 
\end{equation}
\item[3.] Time derivative estimate:
\begin{equation}
\label{TD-AS}
\begin{split}
&\|\partial_t \phi_m \|_{L^\infty(0,T;(H^1(\Omega))')}^2
+\int_0^T \| \nabla \partial_t \phi_m(\tau)\|_{L^2(\Omega)}^2 \, \d \tau\\
&\quad \leq C\Big(1+\| \nabla \mu_0\|_{L^2(\Omega)}^2+ \| \vv\|_{L^\infty(0,T;L^2(\Omega))}^2+ \int_0^T \| \partial_t \vv(\tau)\|_{L^2(\Omega)}^2 \, \d \tau \Big) \mathrm{e}^{C\int_0^T 1+\| \vv(\tau)\|_{L^2(\Omega)}^2\, \d \tau},
\end{split}
\end{equation}
where the constant $C$ only depends on $\Omega$ and $\theta_0$.
\end{itemize}
\smallskip

Next, we look for the approximated velocity field
$$
\uu_m(x,t)= \sum_{j=1}^m a_j^m(t) \ww_j(x)
$$
that solves the Galerkin approximation of \eqref{AGG}$_1$ as follows
\begin{align}
\label{w-ap1}
\begin{split}
&\begin{aligned}[t]
(\rho(\phi_m) \partial_t \uu_m, \ww_l)&+(\rho(\phi_m)(\vv\cdot \nabla)\uu_m,\ww_l)+(\nu(\phi_m) D \uu_m, \nabla \ww_l)\\
&-\frac{\rho_1-\rho_2}{2} ( (\nabla \mu_m \cdot \nabla) \uu_m, \ww_l)= (\nabla \phi_m \otimes \nabla \phi_m, \nabla \ww_l), \quad \forall \, l=1,\dots,m,
\end{aligned}
\end{split}
\end{align}
which is completed with the initial condition $\uu_m(\cdot, 0)= \mathbb{P}_m\uu_0$.
Setting $\mathbf{A}^m(t)=(a_1^m(t), \dots, a_m^m(t))^T$, \eqref{w-ap1} is equivalent to the system of differential equations
\begin{equation}
\label{AS-ODE}
\mathbf{M}^m(t) \ddt \mathbf{A}^m= \mathbf{L}^m(t) \mathbf{A}^m + \mathbf{G}^m(t),
\end{equation}
where the matrices $\mathbf{M}^m(t)$, $\mathbf{L}^m(t)$ and the vector $\mathbf{G}^m(t)$ are given by
\begin{align*}
&(\mathbf{M}^m(t))_{l,j}= \int_{\Omega} \rho(\phi_m) \ww_l \cdot \ww_j \, \d x,\\
&(\mathbf{L}^m(t))_{l,j}=\int_{\Omega} 
\Big( \rho(\phi_m) (\vv \cdot \nabla) \ww_j \cdot \ww_l + \nu(\phi_m) D \ww_j : \nabla \ww_l - \Big(\frac{\rho_1-\rho_2}{2}\Big) (\nabla \mu_m \cdot \nabla) \ww_j \cdot \ww_l \Big) \, \d x,\\
&(\mathbf{G}^m(t))_l= \int_{\Omega} \nabla \phi_m \otimes \nabla \phi_m : \nabla \ww_l \, \d x,
\end{align*}
and the initial condition is 
$$
\mathbf{A}^m(0)=((\mathbb{P}_m\uu_{0}, \ww_1), \dots, (\mathbb{P}_m \uu_{0},\ww_m))^T.
$$
Thanks to \eqref{phi-As}, it follows that $\phi_m \in C([0,T]; W^{1,4}(\Omega))$.  This, in turn, implies that 
$\rho(\phi_m), \nu(\phi) \in C(\overline{\Omega \times [0,T]})$. In addition, we recall that $\vv \in C([0,T]; \H_\sigma)$ and $\nabla \mu_m\in C([0,T];L^2(\Omega))$. As a consequence, it follows that $\mathbf{M}^m$ and $\mathbf{L}^m$ belong to $C([0,T];\mathbb{R}^{m \times m})$, and $\mathbf{G}^m \in C([0,T];\mathbb{R}^m)$. Furthermore, the matrix $\mathbf{M}^m(\cdot)$ is definite positive on $[0,T]$, and so the inverse $(\mathbf{M}^m)^{-1} \in C([0,T]; \mathbb{R}^{m\times m})$. Therefore, the classical existence and uniqueness theorem for system of linear ODEs entails that there exists a unique vector 
 $\mathbf{A}^m \in C^1([0,T];\mathbb{R}^m)$ that solves \eqref{AS-ODE} on $[0,T]$. This implies that the problem \eqref{w-ap1} has a unique solution $\uu_m \in C^1([0,T];\V_m)$.

Next, multiplying \eqref{w-ap1} by $a_l^m$ and summing over $l$, we find
\begin{align*}
\int_{\Omega} \rho(\phi_m) \partial_t \Big( \frac{|\uu_m|^2}{2} \Big) \, \d x&+ \int_{\Omega} \rho(\phi_m) \vv \cdot \nabla \Big(  \frac{|\uu_m|^2}{2} \Big) \, \d x + \int_{\Omega} \nu(\phi_m) |D \uu_m|^2 \, \d x\\
&- \frac{\rho_1-\rho_2}{2} \int_{\Omega} \nabla \mu_m \cdot \nabla \Big(  \frac{|\uu_m|^2}{2} \Big) \, \d x = \int_{\Omega} \nabla \phi_m \otimes \nabla \phi_m : \nabla \uu_m \, \d x. 
\end{align*}
By integration by parts, we have
\begin{align*}
\ddt \int_{\Omega} \rho(\phi_m) \frac{|\uu_m|^2}{2}\, \d x &-
\int_{\Omega} \Big( \partial_t \rho(\phi_m)+ \div \big(\rho(\phi_m) \vv \big) \Big)  \frac{|\uu_m|^2}{2}\, \d x
+ \int_{\Omega} \nu(\phi_m) |D \uu_m|^2 \, \d x\\
&+ \frac{\rho_1-\rho_2}{2} \int_{\Omega}\Delta \mu_m  \frac{|\uu_m|^2}{2} \, \d x = \int_{\Omega} \nabla \phi_m \otimes \nabla \phi_m : \nabla \uu_m \, \d x. 
\end{align*}
Since $\rho'(\phi_m)= \frac{\rho_1-\rho_2}{2}$ and $\div \vv=0$, by using \eqref{cCH}$_1$, we observe that
\begin{align*}
&-\int_{\Omega} \Big( \partial_t \rho(\phi_m)+ \div \big(\rho(\phi_m) \vv \big) \Big)  \frac{|\uu_m|^2}{2}\, \d x+ \frac{\rho_1-\rho_2}{2} \int_{\Omega}\Delta \mu_m  \frac{|\uu_m|^2}{2} \, \d x \\
&\qquad = \int_{\Omega} \rho'(\phi_m) \Big( \partial_t \phi_m + \vv \cdot \nabla \phi_m -\Delta \mu_m\Big)  \frac{|\uu_m|^2}{2}\, \d x=0.
\end{align*}
Thus, we deduce that
\begin{equation}
\ddt \int_{\Omega} \rho(\phi_m) \frac{|\uu_m|^2}{2}\, \d x 
+ \int_{\Omega} \nu(\phi_m) |D \uu_m|^2 \, \d x\\
= \int_{\Omega} \nabla \phi_m \otimes \nabla \phi_m : \nabla \uu_m \, \d x. 
\end{equation}
By using \eqref{GN-L4}, \eqref{KORN} and \eqref{phi-As}, we have 
$$
- \int_{\Omega} \nabla \phi_m \otimes \nabla \phi_m : \nabla \uu_m \, \d x \leq \| \nabla \phi_m \|_{L^4(\Omega)}^2 \| \nabla \uu_m \|_{L^2(\Omega)}\leq \frac{\nu_\ast}{2} \| D \uu_m\|_{L^2(\Omega)}^2 +  C \| \phi_m \|_{H^2(\Omega)}^2,
$$
for some constant $C$ depending only on $\Omega$ and $\nu_\ast$.
Since $\| \phi_m\|_{H^2(\Omega)}\leq C(1+ \| \Delta \phi_m\|_{L^2(\Omega)})$, we arrive at
\begin{equation}
\ddt \int_{\Omega} \rho(\phi_m) \frac{|\uu_m|^2}{2}\, \d x 
+ \frac{\nu_\ast}{2}\int_{\Omega} |D \uu_m|^2 \, \d x\\
\leq C\big(1+ \| \Delta \phi_m\|_{L^2(\Omega)}^2 \big).
\end{equation}
In light of \eqref{L2-AS}, we infer that
$$
\sup_{t\in [0,T]}  \int_{\Omega} \rho(\phi_m(t)) \frac{|\uu_m(t)|^2}{2}\, \d x \leq  \int_{\Omega} \rho(\phi_0) \frac{|\mathbb{P}_m \uu_0|^2}{2}\, \d x + C \Big( T+ \| \phi_0\|_{L^2(\Omega)}^2 + \frac{\theta_0^2}{2}T \Big).
$$
This, in turn, implies that 
\begin{equation}
\label{uL2-AS}
\| \uu_m\|_{C([0,T];\H_\sigma)}\leq R_0,
\end{equation}
where the constant $R_0$ depends on $\rho_\ast$, $\rho^\ast$, $\nu_\ast$, $\theta_0$, $\|\uu_0 \|_{L^2(\Omega)}$, $T$, $\Omega$.
As an immediate consequence, we deduce that
\begin{equation}
\| \uu_m \|_{L^2(0,T;\H_\sigma)} \leq R_0 \sqrt{T}=:R_1.
\end{equation}

Next, we proceed in estimating the time derivative of $\uu_m$. To this aim, multiplying \eqref{w-ap1} by $\ddt a_l^m$ and summing over $l$, we obtain
\begin{align*}
\rho_\ast \| \partial_t \uu_m \|_{L^2(\Omega)}^2 
&\leq - (\rho(\phi_m)(\vv\cdot \nabla)\uu_m,\partial_t \uu_m)-(\nu(\phi_m) D \uu_m, \nabla \partial_t \uu_m)\\
&\quad +\frac{\rho_1-\rho_2}{2} ( (\nabla \mu_m \cdot \nabla) \uu_m, \partial_t \uu_m) +(\nabla \phi_m \otimes \nabla \phi_m, \nabla \partial_t \uu_m).
\end{align*}
By exploiting \eqref{GN-L4} and \eqref{Rev-SI}, we find
\begin{align*}
\rho_\ast \| \partial_t \uu_m \|_{L^2(\Omega)}^2 
&\leq \rho^\ast \| \vv \|_{L^2(\Omega)} \| \nabla\uu_m \|_{L^\infty(\Omega)} \| \partial_t \uu_m\|_{L^2(\Omega)} +
\nu^\ast \|D \uu_m\|_{L^2(\Omega)} \| \nabla \partial_t \uu_m \|_{L^2(\Omega)}\\
&\quad +\Big| \frac{\rho_1-\rho_2}{2} \Big| \| \nabla \uu_m \|_{L^\infty(\Omega)}
\| \nabla \mu_m \|_{L^2(\Omega)} \| \partial_t \uu_m \|_{L^2(\Omega)} + \| \nabla \phi_m \|_{L^4(\Omega)}^2 \|\nabla \partial_t \uu_m\|_{L^2(\Omega)}\\
&\leq   \rho^\ast C  \| \vv\|_{\L2} \| \uu_m\|_{H^3(\Omega)} \| \partial_t \uu_m\|_{L^2(\Omega)} +  \nu^\ast C_m^2 \| \uu_m\|_{\L2} \| \partial_t \uu_m\|_{L^2(\Omega)} \\
&\quad + C \Big| \frac{\rho_1-\rho_2}{2} \Big| \| \uu_m\|_{H^3(\Omega)} \| \nabla \mu_m\|_{\L2}  \| \partial_t \uu_m \|_{L^2(\Omega)}+ C_m \| \phi_m\|_{H^2(\Omega)}  \| \partial_t \uu_m \|_{L^2(\Omega)}\\
&\leq   \rho^\ast C_m  \| \vv\|_{\L2} \| \uu_m\|_{\L2} \| \partial_t \uu_m\|_{\L2} +  \nu^\ast C_m^2 \| \uu_m\|_{\L2} \| \partial_t \uu_m\|_{\L2} \\
&\quad + C_m \Big| \frac{\rho_1-\rho_2}{2} \Big| \| \uu_m\|_{\L2)} \| \nabla \mu_m\|_{\L2}  \| \partial_t \uu_m \|_{\L2}
+ C_m \| \phi_m\|_{H^2(\Omega)}  \| \partial_t \uu_m \|_{\L2}.
\end{align*}
Then, by \eqref{L2-AS}, \eqref{EE-AS}, \eqref{uL2-AS} we eventually infer that
\begin{align*}
&\int_0^T \| \partial_t \uu_m(\tau)\|_{\L2}^2 \, \d \tau \\
&\quad\leq \Big( \frac{\rho^\ast}{\rho_\ast} C_m R_0 \Big)^2 \int_0^T \| \vv(\tau)\|_{\L2}^2 \, \d \tau + \Big( \frac{\nu^\ast}{\rho_\ast} C_m^2 R_0 \Big)^2 T \\
&\qquad + \Big(C_m \Big| \frac{\rho_1-\rho_2}{2} \Big| R_0 \Big)^2 
\int_0^T \| \nabla \mu(\tau)\|_{\L2}^2 \, \d \tau  
+ \Big( \frac{C_m}{\rho_\ast} \Big)^2 \int_0^T (1+ \| \Delta \phi_m(\tau)\|_{\L2}^2) \, \d \tau\\
&\quad\leq \Big[ \Big( \frac{\rho^\ast}{\rho_\ast} C_m R_0 \Big)^2
+ \Big(C_m \Big| \frac{\rho_1-\rho_2}{2} \Big| R_0 \Big)^2 \Big]
 \int_0^T \| \vv(\tau)\|_{\L2}^2 \, \d \tau 
+ \Big( \frac{\nu^\ast}{\rho_\ast} C_m^2 R_0 \Big)^2 T \\
&\qquad + \Big(C_m \Big| \frac{\rho_1-\rho_2}{2} \Big| R_0 \Big)^2 
  \Big( 2  E_{\text{free}}(\phi_0) + \theta_0 | \Omega| + \frac{\theta_0^3}{2}T \Big) 
+C \Big( \frac{C_m}{\rho_\ast} \Big)^2 \Big( 1+ |\Omega| + \frac{\theta_0^2}{2}T \Big).
\end{align*}
Thus, there exist two positive constants $R_2$ and $R_3$, depending only on $\rho_\ast$, $\rho^\ast$, $\nu_\ast$, $\theta_0$, $\|\uu_0 \|_{L^2(\Omega)}$, $E_{\text{free}}(\phi_0)$, $T$, $\Omega$, $m$,  
such that 
\begin{equation}
\label{utL2-AS}
\int_0^T \| \partial_t \uu_m(\tau)\|_{\L2}^2 \, \d \tau
\leq R_2  \int_0^T \| \vv(\tau)\|_{\L2}^2 \, \d \tau+ R_3.
\end{equation}

We are now in a position to state the setting of the fixed point argument. Let us define $R_4= \sqrt{R_2 R_1^2+R_3}$. We introduce the set
$$
S=\Big\lbrace \uu \in W^{1,2}(0,T;\V_m): \| \uu\|_{L^2(0,T;\V_m)}\leq R_1, \ \| \partial_t \uu\|_{L^2(0,T;\V_m)}\leq R_4 \Big\rbrace \subset L^2(0,T;\V_m), 
$$
and the map 
$$
\Lambda: S \rightarrow L^2(0,T;\V_m), \quad \Lambda(\vv)= \uu_m,
$$
where $\uu_m$ is the solution to the system \eqref{w-ap1}.Thanks to \eqref{uL2-AS} and \eqref{utL2-AS}, we deduce that
$$
\Lambda: S \rightarrow S. 
$$
We notice that $S$ is a convex and compact set in $L^2(0,T;\V_m)$. 

We are left to prove that the map $\Lambda$ is continuous. Let us consider a sequence $\lbrace \vv_n\rbrace\subset S$ such that $\vv_n \rightarrow \widetilde{\vv}$ in $L^2(0,T;\V_m)$. By arguing as above, there exists a sequence $\lbrace (\psi_n, \mu_n) \rbrace$ and $(\widetilde{\psi}, \widetilde{\mu})$ that solve the convective Cahn-Hilliard equation \eqref{cCH}-\eqref{cCH-c}, where $\vv$ is replaced by $\vv_n$ and $\widetilde{\vv}$, respectively. Since $\lbrace \vv_n \rbrace $ and $\widetilde{\vv}$ belong to $S$, and $ E_{\text{free}}(\phi_0)<\infty$, we infer from \cite[Theorem 6]{A2009} that 
\begin{equation}
\label{wc1-AS}
\| \psi_n -\widetilde{\psi}\|_{L^\infty(0,T;(H^1(\Omega))')} \rightarrow 0, \quad \text{as } \ n \rightarrow \infty.
\end{equation}
On the other hand, using again that $\lbrace \vv_n \rbrace $ and $\widetilde{\vv}$ belong to $S$, together with the continuous embedding $W^{1,2}(0,T;\V_m) \hookrightarrow C([0,T];\V_m)$, it follows from \eqref{TD-AS} that
\begin{equation}
\label{td-AS}
\begin{split}
&\|\partial_t \psi_n \|_{L^\infty(0,T;(H^1(\Omega))')}
+\| \partial_t \psi_n\|_{L^2(0,T;H^1(\Omega)} \leq C,\\
&\|\partial_t \widetilde{\psi} \|_{L^\infty(0,T;(H^1(\Omega))')}
+\| \partial_t \widetilde{\psi}\|_{L^2(0,T;H^1(\Omega)} \leq C,
\end{split}
\end{equation}
for some constant $C$ which depends on $\phi_0$, $T$, $R_1$, $R_4$, $\theta_0$, $\Omega$, but is independent of $n$.
By comparison in \eqref{cCH}$_1$, it is easily seen that 
$$
\| \mu_n\|_{L^\infty(0,T;H^1(\Omega))}\leq C, \quad 
\| \widetilde{\mu}\|_{L^\infty(0,T;H^1(\Omega))}\leq C.
$$
By exploiting \cite[Lemma A.4 and Lemma A.6]{CG2020}, we obtain
\begin{align*}
&\| \psi_n\|_{L^\infty(0,T;W^{2,p}(\Omega))}+ 
\| F'(\psi_n)\|_{L^\infty(0,T;L^p(\Omega))}
+\| F''(\psi_n)\|_{L^\infty(0,T;L^p(\Omega))}\leq C_p,\\
&\| \widetilde{\psi}\|_{L^\infty(0,T;W^{2,p}(\Omega))}
+\| F'(\widetilde{\psi})\|_{L^\infty(0,T;L^p(\Omega))}
+ \| F''(\widetilde{\psi})\|_{L^\infty(0,T;L^p(\Omega))}
\leq C_p,
\end{align*}
for all $p\in [2,\infty)$, where the constant $C_p$ depends on $p$, 
$\phi_0$, $T$, $R_1$, $R_4$, $\theta_0$, $\Omega$, but is independent of $n$. Thanks to the above estimates, we infer that
$$
\| F'(\psi_n)\|_{L^\infty(0,T;H^1(\Omega))}\leq C, \quad \| F'(\widetilde{\psi})\|_{L^\infty(0,T;H^1(\Omega))}\leq C,
$$
which, in turn, gives us 
\begin{equation}
\label{H3-AS}
\| \psi_n\|_{L^\infty(0,T;H^3(\Omega))} \leq C, \quad
\| \widetilde{\psi}\|_{L^\infty(0,T;H^3(\Omega))}\leq C,
\end{equation}
for some constant $C$ independent of $n$. By standard interpolation, we deduce from \eqref{wc1-AS} and \eqref{H3-AS} that 
\begin{equation}
\label{H2L-AS}
\| \psi_n -\widetilde{\psi}\|_{L^\infty(0,T;H^2(\Omega))} \rightarrow 0, \quad \text{as } \ n \rightarrow \infty.
\end{equation}
As a consequence, by using the definition of $\mu_n-\widetilde{\mu}$ and the above estimates, we eventually obtain
\begin{equation}
\label{MUL-AS}
\| \mu_n -\widetilde{\mu}\|_{L^\infty(0,T;L^2(\Omega))} \rightarrow 0, \quad \text{as } \ n \rightarrow \infty.
\end{equation}
Next, we introduce $\uu_n = \Lambda( \vv_n) \in S$, for any $n \in \mathbb{N}$, and $\widetilde{\uu}=\Lambda(\widetilde{\vv})\in S$. 
We define $ \uu= \uu_n-\widetilde{\uu}$, $ \psi= \psi_n-\widetilde{\psi}$, $ \vv=\vv_n - \widetilde{\vv}$, and $ \mu= \mu_n-\widetilde{\mu}$. We have the system
\begin{align}
\label{w-diff}
\begin{split}
&(\rho(\psi_n) \partial_t \uu, \ww)+ ((\rho(\psi_n)-\rho(\widetilde{\psi})) \partial_t \widetilde{\uu}, \ww)
+(\rho(\psi_n)(\vv_n\cdot \nabla)\uu_n - \rho(\widetilde{\psi})
(\widetilde{\vv}\cdot \nabla) \widetilde{\uu},\ww)\\
&\qquad +(\nu(\psi_n) D  \uu, \nabla \ww)
+ ((\nu(\psi_n)-\nu(\widetilde{\psi})) D \widetilde{\uu}, \nabla \ww)\\
&\qquad-\frac{\rho_1-\rho_2}{2} ( (\nabla \mu_n \cdot \nabla) \uu_n- 
(\nabla \widetilde{\mu}\cdot \nabla) \widetilde{\uu}, \ww)= (\nabla \psi_n \otimes \nabla \psi_n- \nabla \widetilde{\psi} \otimes \nabla \widetilde{\psi}, \nabla \ww), 
\end{split}
\end{align}
for all $\ww \in \V_m$, for all $t \in [0,T]$. Taking $\ww=\uu$, we obtain
\begin{align*}
&\frac12 \ddt \int_{\Omega} \rho(\psi_n) |\uu|^2 \, \d x
+\int_{\Omega} \nu(\psi_n)| D \uu|^2 \, \d x \\
&\quad =
 \frac{\rho_1-\rho_2}{4} \int_{\Omega}  \partial_t \psi_n |\uu|^2 \, \d x-  \frac{\rho_1-\rho_2}{2} \int_{\Omega}  \psi  (\partial_t \widetilde{\uu} \cdot \uu) \, \d x \\
&\qquad - \int_{\Omega} \big(\rho(\psi_n)(\vv_n\cdot \nabla)\uu_n - \rho(\widetilde{\psi})
(\widetilde{\vv}\cdot \nabla) \widetilde{\uu} \big) \cdot \uu \, \d x- \frac{\nu_1-\nu_2}{2} \int_{\Omega}  \psi (D \widetilde{\uu} : D \uu) \, \d x\\
&\qquad+  \frac{\rho_1-\rho_2}{2} \int_{\Omega}
\big( (\nabla \mu_n \cdot \nabla) \uu_n - (\nabla \widetilde{\mu} \cdot \nabla) \widetilde{\uu} \big) \cdot \uu \, \d x  +\int_\Omega \big(\nabla \psi_n \otimes \nabla \psi+ \nabla \psi \otimes \nabla \widetilde{\psi} \big) : \nabla  \uu \, \d x.
\end{align*}
By using \eqref{KORN} and the Sobolev embedding, we have 
\begin{align*}
 \frac{\rho_1-\rho_2}{4}  \int_{\Omega} \partial_t \psi_n |\uu|^2 \, \d x
&\leq C \| \partial_t \psi_n\|_{L^6(\Omega)} \| \uu\|_{\L2} \| \uu \|_{L^3(\Omega)}
\\
& \leq \frac{\nu_\ast}{10} \| D \uu\|_{\L2}^2 + C \| \partial_t \psi_n\|_{H^1(\Omega)}^2 \| \uu\|_{L^2}^2,
\end{align*}
and
\begin{align*}
-\frac{\rho_1-\rho_2}{2} \int_{\Omega} \psi  (\partial_t \widetilde{\uu} \cdot \uu) \, \d x 
&\leq C \| \psi\|_{L^\infty(\Omega)} \| \partial_t \widetilde{\uu}\|_{\L2} \| \uu \|_{\L2}\\
&\leq C\| \uu\|_{L^2(\Omega)}^2+ C \| \psi\|_{H^2(\Omega)}^2.
\end{align*}
Since $\vv_n$, $\widetilde{\vv}$ and $\uu_n $ belong to $S$, by \eqref{KORN} and \eqref{Rev-SI} we get
\begin{align*}
&-\int_{\Omega} \big(\rho(\psi_n)(\vv_n\cdot \nabla)\uu_n - \rho(\widetilde{\psi}) (\widetilde{\vv}\cdot \nabla) \widetilde{\uu} \big) \cdot \uu \, \d x\\
&= - \frac{\rho_1-\rho_2}{2} \int_{\Omega} \psi ((\vv_n\cdot \nabla) \uu_n) \cdot \uu \,\ d x- \int_{\Omega} \rho(\widetilde{\psi})((\vv \cdot \nabla) \uu_n) \cdot \uu \, \d x - \int_{\Omega} \rho(\widetilde{\psi}) ((\widetilde{\vv}\cdot \nabla) \uu) \cdot \uu \, \d x \\
&\leq C \| \psi\|_{L^\infty(\Omega)} \| \vv_n\|_{L^\infty(\Omega)}  \| \nabla \uu_n\|_{\L2} \| \uu\|_{\L2}  + C \| \vv\|_{\L2} \| \nabla \uu_n\|_{L^\infty(\Omega)} \| \uu\|_{\L2}\\
&\quad + C \| \widetilde{\vv}\|_{L^\infty(\Omega)} \| \nabla \uu \|_{\L2} \| \uu \|_{\L2}\\
&\leq C_m \| \psi\|_{H^2(\Omega)} \| \uu\|_{\L2}+ C_m \| \vv\|_{\L2} 
\| \uu\|_{\L2}+ C \| \nabla \uu\|_{\L2} \| \uu\|_{\L2}\\
&\leq   \frac{\nu_\ast}{10} \| D \uu\|_{\L2}^2 + C_m \| \uu\|_{\L2}^2+
+ C_m \| \psi\|_{H^2(\Omega)}^2 + C_m \| \vv\|_{\L2}^2.
\end{align*}
In a similar way, we find
\begin{align*}
- \frac{\nu_1-\nu_2}{2} \int_{\Omega} \psi (D \widetilde{\uu} : D \uu) \, \d x & \leq C \| \psi\|_{L^\infty(\Omega)} \| D \widetilde{\uu}\|_{\L2} \| D \uu \|_{\L2}\\
&\leq  \frac{\nu_\ast}{10} \| D \uu\|_{\L2}^2 + C_m \| \psi\|_{H^2(\Omega)}^2,
\end{align*}
and
\begin{align*}
& \frac{\rho_1-\rho_2}{2} \int_{\Omega}
\big( (\nabla \mu_n \cdot \nabla) \uu_n - (\nabla \widetilde{\mu} \cdot \nabla) \widetilde{\uu} \big) \cdot \uu \, \d x \\
&=- \frac{\rho_1-\rho_2}{2} \int_{\Omega} (\mu_n \Delta \uu_n-\widetilde{\mu} \Delta \widetilde{\uu} ) \cdot \uu \, \d x
-\frac{\rho_1-\rho_2}{2} \int_{\Omega} (\mu_n \nabla \uu_n-\widetilde{\mu} \nabla \widetilde{\uu}): \nabla \uu \, \d x\\
&= -\frac{\rho_1-\rho_2}{2} \int_{\Omega} (\mu \Delta \uu_n+\widetilde{\mu} \Delta \uu ) \cdot \uu \, \d x
-\frac{\rho_1-\rho_2}{2} \int_{\Omega} (\mu \nabla \uu_n-\widetilde{\mu} \nabla \uu): \nabla \uu \, \d x\\
&\leq C \| \mu\|_{\L2} \| \Delta \uu_n\|_{\L2} \| \uu\|_{L^\infty(\Omega)}+
C \| \widetilde{\mu}\|_{L^6(\Omega)} \| \Delta \uu\|_{\L2} \| \uu\|_{L^3(\Omega)} \\
&\quad + C \| \mu\|_{\L2} \| \nabla \uu_n\|_{L^6(\Omega)} \| \nabla \uu\|_{L^3(\Omega)} + C \| \widetilde{\mu}\|_{L^6(\Omega)} \| \nabla \uu\|_{L^6(\Omega)} \| \nabla \uu\|_{L^3(\Omega)}\\
& \leq C_m \| \mu\|_{\L2} \| \nabla \uu\|_{\L2} + C_m \| \nabla \uu \|_{\L2} \| \uu\|_{\L2}\\
&\leq  \frac{\nu_\ast}{10} \| D \uu\|_{\L2}^2 + C_m \| \mu\|_{\L2}^2+ C_m \| \uu \|_{\L2}^2.
\end{align*}
By Sobolev embedding and \eqref{H3-AS}, we have
\begin{align*}
\int_\Omega \big(\nabla \psi_n \otimes \nabla \psi+ \nabla \psi \otimes \nabla \widetilde{\psi} \big) : \nabla \uu \, \d x
&\leq C \big( \| \psi_n\|_{H^2(\Omega)}+ \| \widetilde{\psi}\|_{H^2(\Omega)} \big) \| \psi\|_{H^2(\Omega)} \| \nabla \uu\|_{\L2}\\
&\leq  \frac{\nu_\ast}{10} \| D \uu\|_{\L2}^2 + C \| \psi\|_{H^2(\Omega)}^2.
\end{align*}
Combining the above inequalities, we arrive at the differential inequality
\begin{align*}
& \ddt \int_{\Omega} \rho(\psi_n) |\uu|^2 \, \d x
\leq h_1(t) \int_{\Omega} \rho(\psi_n) |\uu|^2 \, \d x + h_2(t),
\end{align*}
where
$$
h_1(t)= C_m \big(1+ \| \partial_t \psi_n(t)\|_{H^1(\Omega)}^2\big),\quad
h_2(t)= C_m\big(\| \psi(t)\|_{H^2(\Omega)}^2 + \| \vv (t)\|_{\L2}^2+ \| \mu(t)\|_{\L2}^2\big).
$$
Therefore, an application of the Gronwall lemma yields
$$
\sup_{t \in [0,T]} \| \uu(t)\|_{\L2}^2 \leq \frac{1}{\rho_\ast} \mathrm{e}^{\int_0^T h_1(\tau) \d \tau} \int_0^T h_2(\tau) \, \d \tau.
$$
Owing to \eqref{td-AS}, \eqref{H2L-AS}, \eqref{MUL-AS}, and the convergence $\vv_n \rightarrow \widetilde{\vv}$ in $L^2(0,T;\V_m)$, we deduce that $\uu_n \rightarrow \widetilde{\uu}$ in $L^\infty(0,T;\V_m)$, which entails that the map $\Lambda$ is continuous.
Finally, we conclude from the Schauder fixed point theorem that the map $\Lambda$ has a fixed point in $S$. This implies the existence of the approximate solution $(\uu_m, \phi_m)$ on $[0,T]$  satisfying \eqref{reg-AS}-\eqref{bc-AS} for any $m \in \mathbb{N}$.

\subsection{A priori Estimates for the Approximate Solutions}
First, we observe that 
\begin{equation}
\label{cons-mass}
\int_{\Omega} \phi(t) \, \d x= \int_{\Omega} \phi_0 \, \d x, \quad \forall \, t \in [0,T].
\end{equation}
Taking $\ww=\uu_m$ in \eqref{NS-w} and integrating by parts, we obtain
\begin{align*}
&\ddt \int_{\Omega} \frac12 \rho(\phi_m) |\uu_m|^2 \, \d x 
 + \int_{\Omega} \nu(\phi_m)|D \uu_m|^2 \, \d x\\
&\quad= \int_{\Omega} \rho'(\phi_m) \big( \partial_t \phi_m  + \uu_m \cdot \nabla \phi_m- \Delta \mu_m \big) \frac{|\uu_m|^2}{2} \, \d x- \int_{\Omega}
\div( \nabla \phi_m\otimes \nabla \phi_m ) \cdot \uu_m \, \d x
 \end{align*}
Thanks to \eqref{CH-w}$_1$, the first term in the right-hand side in the above equality is zero. We recall that
\begin{align*}
-\div(\nabla \phi_m \otimes \nabla \phi_m)
&=  - \nabla \Big( \frac12 |\nabla \phi_m|^2\Big)- \Delta \phi_m \nabla \phi_m\\
&= \mu_m \nabla \phi_m -  \nabla \Big( \frac12 |\nabla \phi_m|^2\Big) - \nabla \Psi(\phi_m).
\end{align*}
Then, we have
\begin{equation}
\label{NS-u}
\ddt \int_{\Omega} \frac12 \rho(\phi_m) |\uu_m|^2 \, \d x 
 + \int_{\Omega} \nu(\phi_m)|D \uu_m|^2 \, \d x= 
 \int_{\Omega} \mu_m \nabla \phi_m \cdot \uu_m \, \d x. 
\end{equation}
Multiplying \eqref{cCH} by $\mu_m$, integrating over $\Omega$ and using the definition of $\mu_m$, we get
\begin{equation}
\label{CH-mu}
\ddt \int_{\Omega} \frac12 |\nabla \phi_m|^2 + \Psi(\phi_m) \, \d x
+ \int_{\Omega} |\nabla \mu_m|^2 \, \d x + \int_{\Omega} \uu_m \cdot \nabla \phi_m \mu_m \, \d x=0.
\end{equation}
By summing \eqref{NS-u} and \eqref{CH-mu}, we obtain
\begin{equation}
\label{EE-m}
\ddt E(\uu_m, \phi_m)
 + \int_{\Omega} \nu(\phi_m)|D \uu_m|^2 \, \d x
 + \int_{\Omega} |\nabla \mu_m|^2 \, \d x=0.
\end{equation}
Integrating in time, we find
$$
E(\uu_m(t), \phi_m(t))+ \int_0^t \int_{\Omega} \nu(\phi_m)|D \uu_m|^2 \, \d x + \int_0^t \int_{\Omega} |\nabla \mu_m|^2 \, \d x= E(\mathbb{P}_m\uu_0,\phi_0), \quad \forall \, t \in [0,T].
$$
Since 
$$
E(\mathbb{P}_m\uu_0,\phi_0) \leq \frac{\rho^\ast}{2} \| \uu_0\|_{\L2}^2+  E_{\text{free}}(\phi_0),
$$
and recalling that $\phi_m \in L^\infty(\Omega\times (0,T))$ such that 
$|\phi_m(x,t)|<1$ almost everywhere in $\Omega\times(0,T)$, we deduce that 
\begin{align}
\label{uL2H1}
&\| \uu_m\|_{L^\infty(0,T;\H_\sigma)} + \| \uu_m\|_{L^2(0,T;\V_\sigma)}\leq C,\\
\label{phiH1}
&\| \phi_m \|_{L^\infty(0,T;H^1(\Omega))}\leq C, \\
\label{nmuL2}
& \| \nabla \mu_m\|_{L^2(0,T;\L2)}\leq C,
\end{align}
where the positive constant $C$ depends on $\| \uu_0\|_{\L2}$ and $ E_{\text{free}}(\phi_0)$, but is independent of $m$. Multiplying \eqref{cCH} by $-\Delta \phi_m$ and integrating over $\Omega$, we have
$$
\| \Delta \phi_m \|_{\L2}^2 + \int_{\Omega} F''(\phi_m) |\nabla \phi_m|^2 \, \d x= \int_{\Omega} \nabla \mu_m \cdot \nabla \phi_m \, \d x+ \theta_0 \|\nabla \phi_m\|_{\L2}^2.
$$
Since $F''(s)>0$ for $s\in (-1,1)$, by using \eqref{phiH1}, we get
\begin{equation}
\label{phiH2e}
\| \Delta \phi_m \|_{\L2}^2  \leq C ( 1+ \| \nabla \mu_m\|_{\L2}),
\end{equation}
for some $C$ independent of $m$. Then, it follows from \eqref{nmuL2} that
\begin{equation}
\label{phiH2}
\| \phi_m\|_{L^4(0,T;H^2(\Omega))}\leq C.
\end{equation}
We recall the well-known inequality
\begin{equation}
\label{F'-L1}
\int_{\Omega} |F'(\phi_m)| \, \d x \leq C \int_{\Omega} F'(\phi_m) (\phi_m -\overline{\phi_0}) \, \d x+ C,
\end{equation}
where the constant $C$ depends on $\overline{\phi_0}$.
Then, multiplying \eqref{cCH}$_2$ by $\phi_m - \overline{\phi_0}$ (cf. \eqref{cons-mass}), we obtain
$$
\int_{\Omega} |\nabla \phi_m|^2 \, \d x+ 
\int_{\Omega} F'(\phi_m) (\phi_m -\overline{\phi_0}) \, \d x 
= \int_{\Omega} (\mu-\overline{\mu}) \phi_m \, \d x + \theta_0 \int_{\Omega} \phi_m (\phi_m -\overline{\phi_0}) \, \d x.
$$
By the Poincar\'{e} inequality and \eqref{phiH1}, we find
\begin{equation}
\label{F'-L1e}
\int_{\Omega} F'(\phi_m) (\phi_m -\overline{\phi_0}) \, \d x \leq C(1+ \| \nabla \mu_m\|_{\L2}).
\end{equation}
Since $\overline{\mu_m}= \overline{F'(\phi_m)}- \theta_0 \overline{\phi_0}$, by combining \eqref{F'-L1} and \eqref{F'-L1e}, we have 
$$
|\overline{\mu_m}|\leq C(1+ \| \nabla \mu_m\|_{\L2}).
$$ 
Thanks to \eqref{normH1-2}, we are led to 
\begin{equation}
\label{mu-H1e}
\| \mu_m\|_{H^1(\Omega)}\leq C(1+ \| \nabla \mu_m\|_{\L2}),
\end{equation}
 which, in turn, implies
 \begin{equation}
 \label{mu-H1}
 \| \mu_m\|_{L^2(0,T;H^1(\Omega))}\leq C,
 \end{equation}
for some constant $C$ independent of $m$. 
In addition, using the boundary conditions \eqref{bc-AS} and \eqref{uL2H1}, we deduce that
\begin{equation}
\label{phit}
\| \partial_t \phi_m\|_{(H^1(\Omega))'}\leq C (1+ \|\nabla \mu_m \|_{\L2}),
\end{equation}
which entails that
$$
\| \partial_t \phi_m\|_{L^2(0,T;(H^1(\Omega))')}\leq C.
$$
Furthermore, by using \cite[Lemma 2]{A2009} or \cite[Lemma A.4]{CG2020}, we infer that, for all $p\in (2,\infty)$, 
\begin{equation}
\label{phiW2pe}
\| \phi_m\|_{W^{2,p}(\Omega)}+ \| F'(\phi_m)\|_{L^p(\Omega)}
\leq C_p (1+ \| \nabla \mu_m\|_{\L2}).
\end{equation}
As a consequence, it holds
\begin{equation}
\label{phiW2p}
\| \phi_m \|_{L^2(0,T;W^{2,p}(\Omega))}+ \| F'(\phi_m)\|_{L^2(0,T;L^p(\Omega))} \leq C_p.
\end{equation}
\smallskip

Next, taking $\ww=\partial_t \uu_m$ in \eqref{NS-w} we get
\begin{equation}
\label{High1}
\begin{split}
&\frac12 \ddt \int_{\Omega} \nu(\phi_m) |D \uu_m|^2 \, \d x+ \int_{\Omega} \rho(\phi_m) |\partial_t \uu_m|^2 \, \d x\\
&\quad =- \int_{\Omega} \rho(\phi_m)  ((\uu_m \cdot \nabla) \uu_m ) \cdot \partial_t \uu_m \, \d x +\frac{\nu_1-\nu_2}{2} \int_{\Omega}  \partial_t \phi _m |D \uu_m|^2 \, \d x\\
&\qquad + \frac{\rho_1-\rho_2}{2}\int_{\Omega} ((\nabla \mu_m \cdot \nabla) \uu_m ) \cdot \partial_t \uu_m - \int_{\Omega} \Delta \phi_m \nabla \phi_m \cdot \partial_t \uu_m \, \d x.
\end{split}
\end{equation}
Computing the duality between $\partial_t \mu_m$ and \eqref{cCH}, we find
$$
\frac12 \ddt \int_{\Omega} |\nabla \mu_m|^2 \, \d x + \l \partial_t \mu_m, \partial_t \phi_m \r + \l \partial_t \mu_m, \uu_m \cdot \nabla \phi_m\r=0.
$$
Notice that 
\begin{align*}
\l \partial_t \mu_m, \partial_t \phi_m \r= \| \nabla \partial_t \phi_m\|_{\L2}^2+
\int_{\Omega} F''(\phi_m) |\partial_t \phi_m|^2 \,\d x - \theta_0 \| \partial_t \phi_m\|_{\L2}^2
\end{align*}
and
$$
\l \partial_t \mu_m, \uu_m \cdot \nabla \phi_m\r=
\ddt \int_{\Omega} \mu_m \uu_m \cdot \nabla \phi_m\, \d x 
- \int_{\Omega}  \mu_m \partial_t \uu_m \cdot \nabla \phi_m \, \d x
-\int_{\Omega} \mu_m \uu_m \cdot \nabla \partial_t \phi_m \, \d x.
$$
Then, we obtain
\begin{equation}
\label{High2}
\begin{split}
&\ddt \bigg[ \int_{\Omega} \frac12 |\nabla \mu_m|^2 \, \d x + \int_{\Omega} \mu_m \uu_m \cdot \nabla \phi_m \, \d x \bigg] + \| \nabla \partial_t \phi_m\|_{\L2}^2\\
&\quad  \leq \theta_0 \| \partial_t \phi_m \|_{\L2}^2 + \int_{\Omega}  \mu_m \partial_t \uu_m \cdot \nabla \phi_m \, \d x
+ \int_{\Omega} \mu_m \uu_m \cdot \nabla \partial_t \phi_m \, \d x.
\end{split}
\end{equation}
By summing \eqref{High1} and \eqref{High2}, we have
\begin{equation}
\label{High3}
\begin{split}
&\ddt H_m + \rho_\ast \|\partial_t \uu_m\|_{\L2}^2
+ \| \nabla \partial_t \phi_m\|_{\L2}^2\\
&\quad \leq - \int_{\Omega} \rho(\phi_m) ( (\uu_m \cdot \nabla) \uu_m ) \cdot \partial_t \uu_m \, \d x +  \frac{\nu_1-\nu_2}{2} \int_{\Omega}\partial_t \phi _m |D \uu_m|^2 \, \d x\\
&\qquad + \frac{\rho_1-\rho_2}{2} \int_{\Omega} ((\nabla \mu_m \cdot \nabla) \uu_m ) \cdot \partial_t \uu_m - \int_{\Omega} \Delta \phi_m \nabla \phi_m \cdot \partial_t \uu_m \, \d x\\
&\qquad  +\theta_0 \| \partial_t \phi_m \|_{\L2}^2 + \int_{\Omega}  \mu_m \partial_t \uu_m \cdot \nabla \phi_m \, \d x
+ \int_{\Omega} \mu_m \uu_m \cdot \nabla \partial_t \phi_m \, \d x\\
&\quad = \sum_{k=1}^7 I_k,
\end{split}
\end{equation}
where
$$
H_m(t)= \frac12 \int_{\Omega} \nu(\phi_m) |D \uu_m|^2 \, \d x +
\frac12 \int_{\Omega} |\nabla \mu_m|^2 \, d x + \int_{\Omega} \mu_m \uu_m \cdot \nabla \phi_m \, \d x.
$$
By \eqref{LADY}, \eqref{KORN}, \eqref{uL2H1}, \eqref{phiH1}, and \eqref{mu-H1e}, 
\begin{align*}
\int_{\Omega} \mu_m \uu_m \cdot \nabla \phi_m \, \d x
&\leq \| \mu_m\|_{L^4(\Omega)} \| \uu_m\|_{L^4(\Omega)} \| \nabla \phi_m \|_{\L2}\\
&\leq C (1+ \| \nabla \mu_m\|_{\L2}) \|\nabla \uu_m \|_{\L2}^\frac12\\
&\leq \frac14 \|\nabla \mu_m\|_{\L2}^2 +\frac{1}{4} \int_{\Omega}\nu(\phi_m) |D \uu_m|^2 \, \d x +C_0,
\end{align*}
for some $C_0$ independent of $m$. Then, we infer that
\begin{equation}
\label{H-b}
H_m \geq  \frac14 \|\nabla \mu_m\|_{\L2}^2 +\frac{1}{4} \int_{\Omega}\nu(\phi_m) |D \uu_m|^2 \, \d x -C_0.
\end{equation}

We now proceed in estimating the terms $I_i$, $i=1,\dots,7$. 
Let $\varpi_1$ and $\varpi_2$ be two positive constant whose values will be determined later. Exploiting \eqref{LADY}, \eqref{KORN} and \eqref{H-b}, we  have
\begin{equation}
\label{I1}
\begin{split}
|I_1| 
&\leq \rho^\ast \| \uu_m \|_{L^4(\Omega)} \| \nabla \uu_m\|_{L^4(\Omega)} \| \partial_t \uu_m\|_{\L2}\\
& \leq \frac{\rho_\ast}{8} \| \partial_t \uu_m\|_{\L2}^2+
C\| \nabla \uu_m \|_{\L2}^2 \| \A \uu_m\|_{\L2}\\
&\leq  \frac{\rho_\ast}{8} \| \partial_t \uu_m\|_{\L2}^2
+\frac{\varpi_1}{3} \| \A \uu_m\|_{\L2}^2+
C\| \nabla \uu_m \|_{\L2}^4.
\end{split}
\end{equation}
By interpolation of Sobolev spaces and \eqref{normH1-2}, \eqref{LADY}, \eqref{phit}, we obtain
\begin{equation}
\label{I2}
\begin{split}
|I_2| 
&\leq C \| \partial_t \phi_m\|_{L^2(\Omega)} \|D \uu_m  \|_{L^4(\Omega)}^2\\
&\leq C \| \partial_t \phi_m \|_{(H^1(\Omega))'}^\frac12 \| \nabla \partial_t \phi_m\|_{\L2}^\frac12 \| D \uu_m \|_{L^2} \| \A \uu_m\|_{\L2}\\ 
&\leq \frac{1}{4} \| \nabla \partial_t \phi_m\|_{\L2}^2 
+C (1+ \| \nabla \mu_m \|_{\L2}^\frac23) \| D \uu_m \|_{\L2}^\frac43 \| \A \uu_m \|_{\L2}^\frac43\\
&\leq \frac{1}{6} \| \nabla \partial_t \phi_m\|_{\L2}^2 + \frac{\varpi_1}{3} \| \A \uu_m\|_{\L2}^2 + C (1+ \| \nabla \mu_m \|_{\L2}^2) \| D \uu_m \|_{\L2}^4.
\end{split}
\end{equation}
By using \eqref{LADY} and \eqref{mu-H1e}, we get
\begin{equation}
\label{I3}
\begin{split}
|I_3| 
&\leq C \| \nabla \uu_m \|_{L^4(\Omega)} \| \nabla \mu_m\|_{L^4(\Omega)} \| \partial_t \uu_m \|_{\L2}\\
& \leq C \| \nabla \uu_m \|_{\L2}^\frac12 \| \A \uu_m\|_{\L2}^\frac12 \| \nabla \mu_m\|_{\L2}^\frac12 \| \mu_m \|_{H^2(\Omega)}^\frac12 \| \partial_t \uu_m\|_{\L2}\\
&\leq C \| \nabla \uu_m \|_{\L2}^\frac12 \| \A \uu_m\|_{\L2}^\frac12 (1+\| \nabla \mu_m\|_{\L2}^\frac34) \| \mu_m \|_{H^3(\Omega)}^\frac14 \| \partial_t \uu_m\|_{\L2}\\
&\leq \frac{\rho_\ast}{4} \| \partial_t \uu_m\|_{\L2}^2 
+ \frac{\varpi_1}{3} \| \A \uu_m\|_{\L2}^2 + \varpi_2 \| \mu_m\|_{H^3(\Omega)}^2\\
&\quad +
C\| \nabla \uu_m\|_{\L2}^4 (1+\| \nabla \mu_m\|_{\L2}^6).
\end{split}
\end{equation}
Exploiting \eqref{phiH2e}, \eqref{phit} and \eqref{phiW2pe}, we find
\begin{equation}
\label{I4}
\begin{split}
|I_4| 
&\leq \| \Delta \phi_m\|_{L^6(\Omega)} \| \nabla \phi_m \|_{L^3(\Omega)}
\| \partial_t \uu_m\|_{\L2}\\
&\leq \frac{\rho_\ast}{4} \| \partial_t \uu_m\|_{\L2}^2 +
C(1+ \| \nabla \mu_m\|_{\L2}^3),
\end{split}
\end{equation}
and
\begin{equation}
\label{I5}
\begin{split}
|I_5|
&\leq C \| \partial_t \phi_m\|_{(H^1(\Omega))'} \| \nabla \partial_t \phi_m \|_{\L2}\\
&\leq \frac16 \|\nabla \partial_t \phi_m \|_{\L2}^2
 +C (1+ \| \nabla \mu_m\|_{\L2}^2).
\end{split}
\end{equation}
Thanks to \eqref{phiH2e} and \eqref{mu-H1e}, we deduce that 
\begin{equation}
\label{I6}
\begin{split}
|I_6| 
&\leq \| \mu_m\|_{L^6(\Omega)} \| \partial_t \uu_m\|_{\L2} \| \nabla \phi_m\|_{L^3(\Omega)}\\
&\leq \frac{\rho_\ast}{4} \| \partial_t \uu_m\|_{\L2}^2
+C (1+ \| \nabla \mu_m\|_{\L2}^3),
\end{split}
\end{equation}
and
\begin{equation}
\label{I7}
\begin{split}
|I_7|
&\leq \| \mu_m \|_{L^6(\Omega)} \|  \uu_m\|_{L^3(\Omega)} \| \nabla \partial_t \phi_m\|_{\L2}\\
&\leq \frac16 \| \nabla \partial_t \phi_m\|_{\L2}^2 +C \| \nabla \uu_m\|_{\L2}^2 (1+\|\nabla \mu_m\|_{\L2}^2).
\end{split}
\end{equation}
Combining \eqref{High3} with \eqref{H-b} and the above estimates of $I_i$, we arrive at
\begin{equation}
\label{High4}
\begin{split}
&\ddt H_m+ \frac{\rho_\ast}{2} \|\partial_t \uu_m\|_{\L2}^2
+ \frac12 \| \nabla \partial_t \phi_m\|_{\L2}^2\\
&\quad \leq  \varpi_1 \| \A \uu_m\|_{\L2}^2 + \varpi_2 \| \mu_m\|_{H^3(\Omega)}^2 + C (1+ (C_0+H_m)^5),
\end{split}
\end{equation}
where the positive constant $C$ depends on the values of $\varpi_1$ and $\varpi_2$ but is independent of $m$. 
We are left to control the norms $\| \A \uu_m\|_{\L2}$ and $\| \mu_m\|_{H^3(\Omega)}$. To this end, taking $\ww=\A \uu_m$ in \eqref{w-ap1}, we have
\begin{equation}
\label{A-H2}
\begin{split}
-\frac12 (\nu(\phi_m) \Delta \uu_m,  \A \uu_m)&=
-(\rho(\phi_m) \partial_t \uu_m, \A\uu_m)
-(\rho(\phi_m)(\uu_m\cdot \nabla)\uu_m,\A \uu_m)\\
&\quad+\frac{\rho_1-\rho_2}{2} ( (\nabla \mu_m \cdot\nabla) \uu_m, \A \uu_m)-
(\Delta \phi_m \nabla \phi_m,  \A \uu_m)\\
&\quad +\frac{\nu_1-\nu_2}{2}(D \uu_m \nabla \phi_m, \A \uu_m).
\end{split}
\end{equation}
By arguing as in \cite{GMT2019}, there exists $\pi_m \in C([0,T];H^1(\Omega))$ such that $-\Delta \uu_m + \nabla \pi_m= \A \uu_m$ almost everywhere in $\Omega \times (0,T)$ and such that 
\begin{equation}
\label{p-est}
\| \pi_m\|_{L^2(\Omega)}\leq C\| \nabla \uu_m\|_{L^2(\Omega)}^\frac12 \| \A \uu_m\|_{L^2(\Omega)}^\frac12,\quad \| \pi_m\|_{H^1(\Omega)}\leq C \| \A \uu_m\|_{L^2(\Omega)},
\end{equation}
where $C$ is independent of $m$. Thus, we deduce that
\begin{align*}
\frac12 (\nu(\phi_m) \A \uu_m,  \A \uu_m)&=
-(\rho(\phi_m) \partial_t \uu_m, \A\uu_m)
-(\rho(\phi_m)(\uu_m\cdot \nabla)\uu_m,\A \uu_m)\\
&\quad+\frac{\rho_1-\rho_2}{2} ( (\nabla \mu_m \cdot \nabla) \uu_m, \A \uu_m)-
(\Delta \phi_m \nabla \phi_m,  \A \uu_m)\\
&\quad +\frac{\nu_1-\nu_2}{2}(D \uu_m \nabla \phi_m, \A \uu_m)
-\frac{\nu_1-\nu_2}{4} (\pi_m \nabla \phi_m, \A \uu_m)\\
&= \sum_{i=8}^{13} I_i.
\end{align*}
By Young's inequality, we have
\begin{align*}
|I_8|
&\leq \rho^\ast \| \partial_t \uu_m\|_{\L2} \| \A \uu_m\|_{\L2}\\
&\leq \frac{\nu_\ast}{24} \| \A \uu_m\|_{\L2}^2 + \frac{6(\rho^\ast)^2}{\nu_\ast} \| \partial_t \uu_m\|_{\L2}^2.
\end{align*}
By using \eqref{LADY}, \eqref{KORN} and \eqref{uL2H1}, we find
\begin{align*}
|I_9|
&\leq \rho^\ast \| \uu_m\|_{L^4(\Omega)} \| \nabla \uu_m\|_{L^4(\Omega)}
\| \A \uu_m\|_{\L2}\\
&\leq C \| \nabla \uu_m\|_{\L2} \| \A \uu_m \|_{\L2}^\frac32\\
&\leq \frac{\nu_\ast}{24} \| \A \uu_m\|_{\L2}^2 +
C \| D \uu_m\|_{\L2}^4,
\end{align*}
and
\begin{align*}
|I_{10}|
&\leq  C \| \nabla \uu_m\|_{L^4(\Omega)} \|\nabla \mu_m\|_{L^4(\Omega)} 
\| \A \uu_m \|_{\L2}\\
&\leq C \| \nabla \uu_m\|_{\L2}^\frac12 \| \A \uu_m\|_{\L2}^\frac32 \| \nabla \mu_m \|_{\L2}^\frac12 \| \mu_m\|_{H^2(\Omega)}^\frac12\\
&\leq \frac{\nu_\ast}{24} \| \A \uu_m\|_{\L2}^2  + C \| \nabla \uu_m\|_{\L2}^2
(1+\| \nabla \mu_m\|_{\L2}^3)\| \mu_m\|_{H^3(\Omega)}\\
&\leq  \frac{\nu_\ast}{24} \| \A \uu_m\|_{\L2}^2  
+\varpi_2 \| \mu_m\|_{H^3(\Omega)}^2
+ C \| \nabla \uu_m\|_{\L2}^4
(1+\| \nabla \mu_m\|_{\L2}^6).
\end{align*}
In light of \eqref{phiH2e} and \eqref{mu-H1e}, we have
\begin{align*}
|I_{11}|
&\leq C \| \Delta \phi_m\|_{L^6(\Omega)}\| \nabla \phi_m \|_{L^3(\Omega)} \| \A \uu_m\|_{\L2}\\
&\leq \frac{\nu_\ast}{24} \| \A \uu_m\|_{\L2}^2  +
C (1+\| \nabla \mu_m\|_{\L2}^3),
\end{align*}
and
\begin{align*}
|I_{12}|
&\leq C \| D \uu_m\|_{\L2} \| \nabla \phi_m\|_{L^\infty(\Omega)} 
\| \A \uu_m\|_{\L2}\\
&\leq  \frac{\nu_\ast}{24} \| \A \uu_m\|_{\L2}^2 
+C (1+\| \nabla \mu_m\|_{\L2}^2) \| D \uu_m\|_{\L2}^2.
\end{align*}
Owing to \eqref{phiW2pe} and \eqref{p-est}, we obtain
\begin{align*}
|I_{13}|
&\leq 
C \| \pi_m\|_{\L2} \| \nabla \phi_m \|_{L^\infty(\Omega)} \| \A \uu_m \|_{\L2}\\
&\leq C \| D \uu_m\|_{\L2}^\frac12 \| \A \uu_m\|_{\L2}^\frac32 (1+\| \nabla \mu_m\|_{\L2})\\
&\leq  \frac{\nu_\ast}{24} \| \A \uu_m\|_{\L2}^2 
+C \| D \uu_m\|_{\L2}^2 (1+\| \nabla \mu_m\|_{\L2}^4).
\end{align*}
Thus, we are led to 
\begin{equation}
\label{Aum}
\frac{\nu_\ast}{4} \| \A \uu_m\|_{\L2}^2\leq \frac{6(\rho^\ast)^2}{\nu_\ast} \| \partial_t \uu_m\|_{\L2}^2+ \varpi_2 \| \mu_m\|_{H^3(\Omega)}^2
+C (1+ (C_0+H_m)^5).
\end{equation}
Next, taking the gradient of \eqref{CH-w}$_1$, and using \eqref{phiW2pe}, we find
\begin{equation}
\label{H31}
\begin{split}
\|\nabla \Delta \mu_m\|_{\L2}
&\leq \| \nabla \partial_t \phi_m\|_{\L2}+
\| \nabla \uu_m \nabla \phi_m\|_{\L2}+ \| \nabla^2 \phi_m \uu_m\|_{\L2}\\
&\leq \| \nabla \partial_t \phi_m\|_{\L2}+ C \| D \uu_m \|_{\L2} \| \nabla \phi_m\|_{L^\infty(\Omega)} + C \| \nabla ^2 \phi_m\|_{L^3(\Omega)} \| \uu_m\|_{L^6(\Omega)}\\
&\leq \| \nabla \partial_t \phi_m\|_{\L2}+ C \| D \uu_m \|_{\L2} (1+ \| \nabla \mu_m\|_{\L2}).
\end{split}
\end{equation}
Since 
$$
\| \mu_m\|_{H^3(\Omega)}^2 \leq C_S (1+ \| \nabla \mu_m\|_{\L2}^2+ \| \nabla \Delta \mu_m\|_{\L2}^2),
$$
for some positive constant $C_S$ independent of $m$, we infer from \eqref{H31} that
\begin{equation}
\label{muH3m}
\| \mu_m\|_{H^3(\Omega)}^2 \leq 2 C_S \| \nabla \partial_t \phi_m\|_{\L2}^2
+C (1+ (C_0+H_m)).
\end{equation}
Let us now set
$$
\varepsilon_1=\frac{\nu_\ast \rho_\ast}{24 (\rho^\ast)^2},\quad
\varepsilon_2=\frac{1}{8C_S},\quad
\varpi_1= \frac12 \Big( \frac{\nu_\ast^2 \rho_\ast}{96 (\rho^\ast)^2}\Big), \quad 
\varpi_2= \frac{1}{16C_S \big( 1+ \frac{\nu_\ast \rho_\ast}{24 (\rho^\ast)^2}  \big)},\quad C_1=1+C_0.
$$
Multiplying \eqref{Aum} and \eqref{muH3m} by 
$\varepsilon_1$ and $\varepsilon_2$, respectively, and summing the resulting inequalities to \eqref{High4}, we deduce the differential inequality
\begin{equation}
\label{High5}
\ddt H_m
+F_m
\leq C (C_1+H_m)^5
\end{equation}
where
$$
F_m(t)= \frac{\rho_\ast}{2} \|\partial_t \uu_m(t)\|_{\L2}^2
+ \frac12 \| \nabla \partial_t \phi_m(t)\|_{\L2}^2
+ \varpi_1 \| \A \uu_m(t)\|_{\L2}^2
+ \varpi_2 \| \mu_m (t)\|_{H^3(\Omega)}^2.
$$
Hence, whenever $\widetilde{T}>0$ satisfies
$$
1-4C\widetilde{T}(C_1+H_m(0))^4>0,
$$
we have
\begin{equation}
H_m(t)\leq \frac{1+H_m(0)}{\big( 1-4C t (C_1+H_m(0))^4  \big)^\frac14}, \quad \forall \, t \in [0,\widetilde{T}].
\end{equation}
We observe that
$$
H_m(0)\leq C_2(\| \uu_0\|_{\V_\sigma}+ \| \mu_0\|_{H^1(\Omega)}),
$$
for a positive constant $C_2$ independent of $m$. Therefore, setting 
$$
T_0=\frac{1}{8C(C_1+C_2(\| \uu_0\|_{\V_\sigma}+ \| \mu_0\|_{H^1(\Omega)}) )^4},
$$
it yields  that
$$
H_m(t)\leq \frac{1+H_m(0)}{2^\frac14}, \quad \forall \, t \in [0,T_0].
$$
Notice that $T_0$ is independent of $m$. 
Thanks to \eqref{H-b}, we infer that
\begin{equation}
\label{HE1}
\sup_{t\in [0,T_0]} \| \nabla \uu_m(t)\|_{\L2} + \sup_{t\in[0,T_0]} \| \nabla \mu_m(t)\|_{\L2}\leq K_1,
\end{equation}
where $K_1$ is a positive constant that depends on $E(\uu_0,\phi_0)$, $\| \uu_0\|_{\V_\sigma}$, $\| \mu_0\|_{H^1(\Omega)}$ and the parameters of the system, but is independent of $m$.
Recalling \eqref{phiW2pe}, and using \cite[Lemma A.6]{CG2020}, we immediately obtain for all $p \in [2,\infty)$
\begin{equation}
\label{HE2}
\sup_{t\in [0,T_0]} \| \phi_m(t)\|_{W^{2,p}(\Omega)} 
+ \sup_{t\in[0,T_0]} \| F'(\phi_m(t))\|_{L^p(\Omega)}
+\sup_{t\in[0,T_0]} \| F''(\phi_m(t))\|_{L^p(\Omega)}
\leq K_{2,p}.
\end{equation}
As a consequence, we have 
\begin{equation}
\label{HE3}
\sup_{t\in [0,T_0]} \| \phi_m(t)\|_{H^3(\Omega)}+ \sup_{t\in[0,T_0]} \| F'''(\phi_m(t))\|_{L^p(\Omega)}
\leq K_3.
\end{equation}
Integrating \eqref{High5} we deduce that
\begin{equation}
\label{HE4}
\int_0^{T_0} 
\|\partial_t \uu_m(\tau)\|_{\L2}^2
+\| \nabla \partial_t \phi_m(\tau)\|_{\L2}^2
+\| \A \uu_m(\tau)\|_{\L2}^2
+\| \mu_m (\tau)\|_{H^3(\Omega)}^2 \, \d \tau \leq K_4.
\end{equation}
Finally, it follows from \eqref{HE2} and \eqref{HE4} that
\begin{equation}
\label{HE5}
\int_0^{T_0} \| \partial_t \mu_m (\tau)\|_{(H^1(\Omega))'}^2 \, \d \tau\leq K_5.
\end{equation}
Here, the constants $K_2$, \dots, $K_5$ depend on the same factors as $K_1$.

\subsection{Passage to the Limit and Existence of Strong Solutions}
We are in a position to pass to the limit as $m\rightarrow \infty$.
More precisely, thanks to the above estimates \eqref{HE1}-\eqref{HE5}, we deduce the following convergences (up to a subsequence)
\begin{equation}
\label{wl-SS}
\begin{split}
\begin{aligned}
&\uu_m \rightharpoonup \uu \quad &&\text{weak-star in } L^\infty(0,T_0;\V_\sigma),\\
&\uu_m \rightharpoonup \uu \quad  &&\text{weakly in } L^2(0,T_0;H^2)\cap W^{1,2}(0,T_0;\H_\sigma),\\
&\phi_m \rightharpoonup \phi \quad  &&\text{weak-star in } L^\infty(0,T_0;H^3(\Omega)),\\
&\phi_m \rightharpoonup \phi \quad  &&\text{weakly in } 
W^{1,2}(0,T_0;H^1(\Omega)),\\
&\mu_m  \rightharpoonup \mu \quad  &&\text{weak-star in } L^\infty(0,T_0;H^1(\Omega)),\\
&\mu_m  \rightharpoonup \mu \quad  &&\text{weakly in } L^2(0,T_0;H^3(\Omega)) \cap W^{1,2}(0,T_0;(H^1(\Omega))').
\end{aligned}
\end{split}
\end{equation}
The strong convergences of $\uu_m$, $\phi_m$ and $\mu_m$ are recovered through the Aubin-Lions lemma. We have
\begin{equation}
\label{sl-SS}
\begin{split}
\begin{aligned}
&\uu_m \rightarrow \uu \quad &&\text{strongly in } L^2(0,T_0;\V_\sigma),\\
&\phi_m \rightarrow \phi \quad  &&\text{strongly in } C([0,T_0];W^{2,p}(\Omega)), \ \forall \, p \in [2,\infty),\\
&\mu_m  \rightarrow \mu \quad  &&\text{strongly in } L^2(0,T_0;H^2(\Omega)).
\end{aligned}
\end{split}
\end{equation}
As a consequence, we infer that
\begin{equation}
\label{sl-SS2}
\begin{split}
\begin{aligned}
&\rho(\phi_m) \rightarrow \rho(\phi) \quad  &&\text{strongly in } C([0,T_0];H^2(\Omega)),\\
&\nu(\phi_m) \rightarrow \nu(\phi) \quad  &&\text{strongly in } C([0,T_0];H^2(\Omega)),\\
&F'(\phi_m) \rightharpoonup F'(\phi) \quad  &&\text{weak-star in } L^\infty(0,T_0;L^p(\Omega)), \ \forall \, p \in [2,\infty),\\
&F''(\phi_m) \rightharpoonup F''(\phi) \quad  &&\text{weak-star in } L^\infty(0,T_0;L^p(\Omega)), \ \forall \, p \in [2,\infty),\\
&F'''(\phi_m) \rightharpoonup F'''(\phi) \quad  &&\text{weak-star in } L^\infty(0,T_0;L^p(\Omega)), \ \forall \, p \in [2,\infty).
\end{aligned}
\end{split}
\end{equation}
The above properties entail the convergence of the nonlinear terms in \eqref{NS-w}. Then, in a standard manner, we pass to the limit as $m\rightarrow \infty$ in \eqref{NS-w}-\eqref{CH-w}. Lastly, since
\begin{align*}
\big( -\rho(\phi) \partial_t \uu - \rho(\phi)(\uu\cdot \nabla)\uu+ \div( \nu(\phi) D \uu)+\rho'(\phi) (\nabla \mu \cdot \nabla) \uu - 
\div(\nabla \phi \otimes \nabla \phi), \ww \big)=0, 
\end{align*}
for all $\ww \in \H_\sigma$,
there exists $P \in L^2(0,T_0;H^1(\Omega))$, $\overline{P}(t)=0$ (see, e.g., \cite{Galdi}) such that
$$
\nabla P=-\rho(\phi) \partial_t \uu - \rho(\phi)(\uu\cdot \nabla)\uu+ \div( \nu(\phi) D \uu)+\rho'(\phi) \nabla \uu \nabla \mu - 
\div(\nabla \phi \otimes \nabla \phi).
$$

\begin{remark}
The proof of Theorem \ref{mr1} holds true in the boundary periodic setting. In particular, the orthogonal dense set in $\HH_\sigma$ can be chosen as the eigenfunctions of the Stokes operator (see \cite{Temam}) augmented by the constant function. 
Moreover, in order to recover the norm of $\uu_m$ in $H^2(\T^2)$ (cf. \eqref{A-H2}), it is sufficient to take $-\Delta \uu_m$ in $\eqref{w-ap1}$ (instead of $\A \uu_m$). In turn, the term $I_{13}$ involving the pressure $\pi_m$ does not appear. The rest of the proof remains valid with few minor changes. 
\end{remark}

\section{Proof of Theorem \ref{mr2}: Global Existence in the Space Periodic Setting}
\label{Glo-Ex}
\setcounter{equation}{0}

In this section we address the global existence of the strong solutions to the AGG system \eqref{AGG} in $\T^2$. 

We consider a strong solution $(\uu,P, \phi)$ to system \eqref{AGG} defined on the maximal interval of existence $(0,T_\ast)$.This satisfies for all $0<T<T_\ast$
\begin{align}
\label{RSS-T}
\begin{split}
&\uu \in C([0,T]; \VV_\sigma) \cap L^2(0,T;\WW_\sigma)\cap W^{1,2}(0,T;\HH_\sigma),\\
&P \in L^2(0,T;H^1(\T^2)),\\
&\phi \in L^\infty(0,T;H^3(\T^2)), \
\partial_t \phi \in L^\infty(0,T;(H^1(\T^2))')\cap L^2(0,T;H^1(\T^2)),\\
&\phi \in L^\infty(\Omega\times (0,T)) : |\phi(x,t)|<1 \ \text{a.e. in } \  \T^2\times(0,T),\\
&\mu \in C([0,T);H^1(\T^2))\cap L^2(0,T;H^3(\T^2))\cap W^{1,2}(0,T;(H^1(\T^2))'), \\
&F'(\phi), F''(\phi), F'''(\phi) \in L^\infty(0,T;L^p(\T^2)),
\end{split}
\end{align}
for all $p\in [2,\infty)$, and 
\begin{equation}
 \label{AGG-T}
\begin{cases}
\rho(\phi)\partial_t \uu + \rho(\phi)(\uu \cdot \nabla) \uu - \rho'(\phi)(\nabla \mu\cdot \nabla) \uu
- \div (\nu(\phi)D\uu) + \nabla P= - \div(\nabla \phi \otimes \nabla \phi)\\
\div \uu=0\\
\partial_t \phi +\uu\cdot \nabla \phi = \Delta \mu\\
\mu= -\Delta \phi+\Psi'(\phi)
\end{cases}
\end{equation}
almost everywhere in $\T^2\times (0,T^\ast)$.
\smallskip

The aim is to show that $T_\ast=\infty$. We assume by contradiction that $T_\ast<\infty$. In the rest of this section, we prove that the norms related to the functional spaces in \eqref{RSS-T} are uniformly bounded on $(0,T_\ast)$. In turn, this entails that $\uu(T_\ast) \in \VV_\sigma$, $\phi(T_\ast) \in H^2(\T^2)$ such that $\| \phi(T_\ast)\|_{L^\infty(\T^2)}\leq 1$, $|\overline{\phi}(T_\ast)|<1$ and $\mu (T_\ast)=-\Delta \phi(T_\ast)+\Psi'(\phi(T_\ast)) \in H^1(\T^2)$. Thus, by the local existence result in Theorem \ref{mr1}, it is possible to extend the solution beyond the time $T_\ast$. As a consequence, the solution exists globally in time.

\subsection{Energy Estimates} We report some basic energy estimates similar to those obtained in Section \ref{Loc-Ex} (cf. \eqref{cons-mass}-\eqref{phiW2p}). 
\smallskip

First, combining \eqref{AGG-T}$_1$ and \eqref{AGG-T}$_3$, the solution satisfies \eqref{AGG}$_1$ almost everywhere in $\T^2\times (0,T^\ast)$. Integrating over $\T^2 \times (0,t)$ with $t<T_\ast$, 
we obtain
\begin{equation}
\label{CMu}
\int_{\T^2} \rho(\phi(t)) \uu(t) \, \d x= \int_{\T^2} \rho(\phi_0) \uu_0 \, \d x, \quad \forall \, t \in [0,T_\ast).
\end{equation}
Similarly, integrating \eqref{AGG-T}$_3$ over $\T^2 \times (0,t)$ with $t<T_\ast$, we get
\begin{equation}
\label{CMphi}
\int_{\T^2} \phi(t) \, \d x= \int_{\T^2} \phi_0 \, \d x, \quad \forall \, t \in [0,T_\ast).
\end{equation}
Thanks to the energy identity \eqref{EE}, we have
$$
E(\uu(T), \phi(T))+ \int_0^T \int_{\T^2} \nu(\phi)|D \uu|^2 \, \d x + \int_0^T \int_{\T^2} |\nabla \mu|^2 \, \d x= E(\uu_0,\phi_0), \quad \forall \, 0\leq T<T_\ast.
$$
Since $E(\uu_0,\phi_0) <\infty$, we find for all $0<T<T_\ast$
\begin{align}
\label{uL2H1-T}
&\| \uu\|_{L^\infty(0,T;L^2(\T^2))} 
+ \| \uu\|_{L^2(0,T;H^1(\T^2))}\leq C,\\
\label{phiH1-T}
&\| \phi\|_{L^\infty(0,T;H^1(\T^2))}\leq C, \\
\label{nmuL2-T}
& \| \nabla \mu\|_{L^2(0,T;\LT)}\leq C.
\end{align}
Here the constant $C$depends on $E_(\uu_0,\phi_0)$, but it is independent of $T_\ast$.
Arguing as in Section \ref{Loc-Ex}, we have
\begin{equation}
\label{phiH2eT}
\|\phi \|_{H^2(\T^2)}  \leq C ( 1+ \| \nabla \mu\|_{\LT}^\frac12),
\end{equation}
and 
\begin{equation}
\label{mu-mT}
|\overline{\mu}|\leq C(1+ \| \nabla \mu\|_{\LT}).
\end{equation}
The latter implies
\begin{equation}
\label{mu-H1eT}
\| \mu\|_{H^1(\T^2)}\leq C(1+ \| \nabla \mu\|_{\LT}).
\end{equation}
In addition, we recall that
\begin{equation}
\label{phitT}
\| \partial_t \phi\|_{(H^1(\T^2))'}\leq C (1+ \|\nabla \mu \|_{\LT}),
\end{equation}
and
\begin{equation}
\label{phiW2peT}
\| \phi\|_{W^{2,p}(\T^2)}+ \| F'(\phi)\|_{L^p(\T^2)}
\leq C_p (1+ \| \nabla \mu\|_{\LT}), \quad \forall \, p \in (2,\infty).
\end{equation}
As a consequence, it follows that, for all $T<T_\ast$,
\begin{align}
\label{phi-T}
&\| \phi\|_{L^4(0,T;H^2(\T^2))}\leq C(1+T), \quad 
\| \phi \|_{L^2(0,T;W^{2,p}(\T^2))}+ \| F'(\phi)\|_{L^2(0,T;L^p(\T^2))} \leq C_p(1+T),\\
\label{mu-phit-T}
&\| \mu\|_{L^2(0,T;H^1(\T^2))}\leq C(1+T), \quad 
\| \partial_t \phi\|_{L^2(0,T;(H^1(\T^2))')}\leq C(1+T).
\end{align}
\smallskip

\subsection{High-Order Estimates for the Concentration}
Taking the duality between $\partial_t \mu$ and \eqref{AGG-T}$_3$, we obtain (cf. \eqref{High2})
\begin{equation}
\label{HO1}
\begin{split}
&\ddt \bigg[ \int_{\T^2} \frac12 |\nabla \mu|^2 \, \d x 
+ \int_{\T^2} \mu \uu \cdot \nabla \phi \, \d x \bigg] 
+ \| \nabla \partial_t \phi\|_{\LT}^2\\
&\quad  \leq \theta_0 \| \partial_t \phi \|_{\LT}^2 
+ \int_{\T^2}  \mu \partial_t \uu \cdot \nabla \phi \, \d x
+ \int_{\T^2} \mu \uu \cdot \nabla \partial_t \phi \, \d x.
\end{split}
\end{equation}
Since
$$
\| \mu\|_{H^3(\T^2)}^2 \leq C_S (1+ \| \nabla \mu\|_{\LT}^2
+ \| \nabla \Delta \mu\|_{\LT}^2),
$$
arguing as in \eqref{H31}, we infer that
\begin{equation}
\label{muH3T}
\| \mu\|_{H^3(\T^2)}^2\leq 2C_S \| \nabla \partial_t \phi\|_{\LT}^2
+C(1+\| \uu\|_{H^1(\T^2)}^2) (1+ \| \nabla \mu\|_{\LT}^2).
\end{equation}
Let us set $\varepsilon=\frac{1}{4C_S}$. Multiplying \eqref{muH3T} by $\varepsilon$ and adding the resulting inequality to \eqref{HO1}, we get
\begin{equation}
\label{HO2}
\begin{split}
&\ddt \bigg[ \int_{\T^2} \frac12 |\nabla \mu|^2 \, \d x 
+ \int_{\T^2} \mu (\uu \cdot \nabla \phi)\, \d x \bigg] 
+ \frac12 \| \nabla \partial_t \phi\|_{\LT}^2+ \varepsilon \| \mu\|_{H^3(\T^2)}^2\\
&\quad  \leq \theta_0 \| \partial_t \phi \|_{\LT}^2 
+ \int_{\T^2}  \mu \partial_t \uu \cdot \nabla \phi \, \d x
+ \int_{\T^2} \mu \uu \cdot \nabla \partial_t \phi \, \d x\\
&\qquad+\varepsilon C (1+\| \uu\|_{H^1(\T^2)}^2) (1+ \| \nabla \mu\|_{\LT}^2).
\end{split}
\end{equation}
By interpolation of Sobolev spaces and \eqref{phitT}
$$
\theta_0 \| \partial_t \phi\|_{\LT}^2\leq \frac18 \| \nabla \partial_t \phi\|_{\LT}^2+ C (1+ \| \nabla \mu\|_{\LT}^2).
$$
By using \eqref{mu-H1eT}, we have
\begin{align*}
\int_{\T^2} \mu \uu \cdot \nabla \partial_t \phi \, \d x
& \leq \| \mu \|_{L^6(\T^2)} \| \uu \|_{L^3(\T^2)} \| \nabla \partial_t \phi\|_{\LT}\\
&\leq \frac18 \| \nabla \partial_t \phi\|_{\LT}^2+ C\| \uu\|_{H^1(\T^2)}^2
(1+ \| \nabla \mu\|_{\LT}^2).
\end{align*}
Thus, we preliminary obtain 
\begin{equation}
\label{HO3}
\begin{split}
&\ddt \bigg[ \int_{\T^2} \frac12 |\nabla \mu|^2 \, \d x 
+ \int_{\T^2} \mu (\uu \cdot \nabla \phi)\, \d x \bigg] 
+ \frac14 \| \nabla \partial_t \phi\|_{\LT}^2+ \varepsilon \| \mu\|_{H^3(\T^2)}^2\\
&\quad \leq \int_{\T^2}  \mu \partial_t \uu \cdot \nabla \phi \, \d x
+C (1+\| \uu\|_{H^1(\T^2)}^2) (1+ \| \nabla \mu\|_{\LT}^2).
\end{split}
\end{equation}
We observe that 
\begin{equation}
\label{Key-term}
\begin{split}
\int_{\T^2}  \mu \partial_t \uu \cdot \nabla \phi \, \d x
&= \int_{\T^2}  \partial_t \uu \cdot (\phi \nabla \mu) \, \d x\\
&=\int_{\T^2}  \rho(\phi)\partial_t \uu \cdot \frac{\phi \nabla \mu}{\rho(\phi)} \, \d x\\	
&=- \int_{\T^2} ((\uu\cdot \nabla) \uu) \cdot \phi \nabla \mu \, \d x
+\int_{\T^2} \rho'(\phi) ((\nabla \mu \cdot \nabla) \uu) \cdot \frac{\phi \nabla \mu}{\rho(\phi)} \, \d x\\	
&\quad + \int_{\T^2} \div(\nu(\phi)D \uu) \cdot \frac{\phi \nabla \mu}{\rho(\phi)} \, \d x -
\int_{\T^2} \nabla P \cdot \frac{\phi \nabla \mu}{\rho(\phi)} \, \d x \\
&\quad - \int_{\T^2} \div(\nabla \phi\otimes \nabla \phi) \cdot \frac{\phi \nabla \mu}{\rho(\phi)} \, \d x\\
&= \int_{\T^2} \uu\otimes \uu : \nabla(\phi \nabla \mu) \, \d x
+\int_{\T^2} \rho'(\phi) ((\nabla \mu \cdot \nabla) \uu) \cdot \frac{\phi \nabla \mu}{\rho(\phi)} \, \d x\\
&\quad -\int_{\T^2} \nu(\phi)D \uu : \nabla\Big( \frac{\phi \nabla \mu}{\rho(\phi)} \Big) \, \d x +
\int_{\T^2} P  \, \div\Big( \frac{\phi \nabla \mu}{\rho(\phi)} \Big) \, \d x \\
&\quad - \int_{\T^2} \div(\nabla \phi\otimes \nabla \phi) \cdot \frac{\phi \nabla \mu}{\rho(\phi)} \, \d x\\
&=\sum_{i=1}^5 W_i.
\end{split}
\end{equation}
Here the periodic boundary conditions played a crucial role to avoid any boundary term. 
We now proceed in estimating the terms $W_i$, $i=1,\dots,5$.
By using \eqref{LADY}, \eqref{uL2H1-T} and \eqref{mu-H1eT}, we have
\begin{equation}
\label{W1}
\begin{split}
|W_1| &\leq C \| \uu\|_{L^4(\T^2)}^2 \big( \| \nabla \phi\|_{L^\infty(\T^2)} 
\| \nabla \mu\|_{\LT}+ 
\|\phi \|_{L^\infty(\T^2)}\| \mu\|_{H^2(\T^2)}\big) \\
&\leq C\| \uu\|_{\LT} \|\uu \|_{H^1(\T^2)} \big( \| \nabla \phi\|_{L^\infty(\T^2)} \| \nabla \mu\|_{\LT} + (1+\| \nabla \mu\|_{\LT}^\frac12) \| \mu\|_{H^3(\T^2)}^\frac12 \big) \\
&\leq \frac{\varepsilon}{8} \| \mu\|_{H^3(\T^2)}^2 
+C (1+\|\uu \|_{H^1(\T^2)}^2+ \| \nabla \phi\|_{L^\infty(\T^2)}^2) (1+\| \nabla \mu\|_{\LT}^2)
\end{split}
\end{equation}
and
\begin{equation}
\label{W2}
\begin{split}
|W_2|
&\leq \Big\| \frac{\phi \rho'(\phi)}{\rho(\phi)} \Big\|_{L^\infty(\T^2)} 
\| \nabla \uu\|_{\LT} \| \nabla \mu\|_{L^4(\T^2)}^2\\
&\leq C \| \nabla \uu\|_{\LT} \| \nabla \mu\|_{\LT} \| \mu\|_{H^2(\T^2)}\\
&\leq C \| \nabla \uu\|_{\LT} (1+\| \nabla \mu\|_{\LT}^\frac32) \| \mu_m\|_{H^3(\T^2)}^\frac12\\
&\leq \frac{\varepsilon}{8} \| \mu\|_{H^3(\T^2)}^2 
+C \|\uu \|_{H^1(\T^2)}^\frac43(1+\| \nabla \mu\|_{\LT}^2).
\end{split}
\end{equation}
We observe that $\rho(\phi)-\phi \rho'(\phi)= \frac{\rho_1+\rho_2}{2}$. 
By interpolation of Sobolev spaces and \eqref{mu-H1eT}, we find
\begin{equation}
\label{W3}
\begin{split}
|W_3| 
&\leq \nu^\ast \| D \uu \|_{\LT} \Big( \Big\|  \frac{\rho(\phi)-\phi \rho'(\phi)}{\rho(\phi)^2} \nabla \mu \otimes \nabla \phi \Big\|_{\LT} + \Big\| \frac{\phi}{\rho(\phi)} \nabla^2 \mu\Big\|_{\LT}\Big)\\
&\leq C \| D \uu\|_{\LT} \| \nabla \mu\|_{\LT} \| \nabla \phi \|_{L^\infty(\T^2)}
+C \| D \uu\|_{\LT} \| \mu\|_{H^2(\T^2)}\\
& \leq C \| D \uu\|_{\LT} \| \nabla \mu\|_{\LT} \| \nabla \phi \|_{L^\infty(\T^2)}
+C \| D \uu\|_{\LT}  (1+\|\nabla \mu \|_{\LT}^\frac12) \| \mu\|_{H^3(\T^2)}^\frac12 \\
&\leq  \frac{\varepsilon}{8} \| \mu\|_{H^3(\T^2)}^2 +
C (1+ \|\uu \|_{H^1(\T^2)}^2+ \| \nabla \phi\|_{L^\infty(\T^2)}^2) (1+\| \nabla \mu\|_{\LT}^2).
\end{split}
\end{equation}
By \eqref{phiH2eT} we obtain
\begin{equation}
\label{W5}
\begin{split}
|W_5|
& \leq C \| \phi\|_{H^2(\T^2)} \| \nabla \phi\|_{L^\infty(\T^2)} \Big\| \frac{\phi}{\rho(\phi)}\Big\|_{L^\infty(\T^2)} \| \nabla \mu\|_{\LT}\\
&\leq C \| \nabla \phi\|_{L^\infty(\T^2)} (1+\|\nabla \mu \|_{\LT}^\frac32 )\\
&\leq C\| \nabla \phi\|_{L^\infty(\T^2)} (1+\|\nabla \mu \|_{\LT}^2).
\end{split}
\end{equation}
Similarly, by \eqref{LADY} and \eqref{mu-H1eT}, we find
\begin{equation}
\label{W4}
\begin{split}
|W_4| 
&\leq \| P\|_{\LT} \Big( \Big\| \frac{\rho(\phi)-\phi\rho'(\phi)}{\rho(\phi)^2} \nabla \phi \cdot \nabla \mu \Big\|_{\LT}+ \Big\|  \frac{\phi}{\rho(\phi)} \Delta \mu  \Big\|_{\LT} \Big)\\
&\leq C \| P\|_{\LT} \big( \| \nabla \phi\|_{L^4(\T^2)} \| \nabla \mu\|_{L^4(\T^2)} + \| \mu\|_{H^2(\T^2)} \big)\\
&\leq C\| P\|_{\LT} \big( \| \nabla \phi\|_{\LT}^\frac12  \| \phi\|_{H^2(\T^2)}^\frac12 \| \nabla \mu\|_{L^2(\T^2)}^\frac12 \| \mu\|_{H^2(\T^2)}^\frac12 + \| \mu\|_{H^2(\T^2)} \big)\\
&\leq C\| P\|_{\LT} \big( \| \phi\|_{H^2(\T^2)}^\frac12 (1+\| \nabla \mu\|_{L^2(\T^2)}^\frac34) \| \mu\|_{H^3(\T^2)}^\frac14 + (1+\| \nabla \mu\|_{\LT}^\frac12) \| \mu\|_{H^3(\T^2)}^\frac12 \big).
\end{split}
\end{equation}
We are now left to find an estimate of the pressure $P$. 
We introduce the function $q$ as the solution to 
\begin{equation}
\label{q-def}
-\div\Big( \frac{\nabla q}{\rho(\phi)} \Big)= P \quad \text{in} \  \T^2 \times (0,T_\ast).
\end{equation}
Since $P\in L^2(0,T; H^1(\T^2))$, for all $0<T<T_\ast$, such that $\overline{P}(t)=0$ for all $t\in (0,T_\ast)$, and $\rho(\phi)\geq \rho_\ast$, the existence of $q$ follows from the Lax-Milgram theorem. In particular, we have $	q \in L^2(0,T_\ast;H^1(\T^2))\cap L^1(0,T_\ast;H^2(\T^2))$, and $\overline{q}(t)=0$ for all $t\in (0,T_\ast)$. In addition, we have the following estimates (cf. \cite[Theorem 2.1]{G2020})
\begin{equation}
\label{q-est}
\| q\|_{H^1(\T^2)}\leq C \| P\|_{\LT}, \quad \| q\|_{H^2(\T^2)}\leq C(1+\| \nabla \phi\|_{L^\infty(\T^2)}) \| P\|_{\LT}.
\end{equation} 
Multiplying \eqref{AGG-T} by $\frac{\nabla q}{\rho(\phi)}$, we find
\begin{align*}
&\int_{\T^2} \div(\uu\otimes \uu) \cdot \nabla q \, \d x
-\int_{\T^2} \frac{\rho'(\phi)}{\rho(\phi)} ((\nabla \mu \cdot \nabla) \uu) \cdot \nabla q \, \d x- \int_{\T^2} \div(\nu(\phi) D \uu) \cdot \frac{\nabla q}{\rho(\phi)} \, \d x\\
&\quad + \int_{\T^2} \nabla P \cdot \frac{\nabla q}{\rho(\phi)} \, \d x=
- \int_{\T^2} \div(\nabla \phi\otimes \nabla \phi) \cdot \frac{\nabla q}{\rho(\phi)} \, \d x.
\end{align*} 
Integrating by parts and using the periodic boundary conditions, and then exploiting \eqref{q-def}, we deduce that
\begin{equation}
\begin{split}
\| P\|_{\LT}^2 
&= \int_{\T^2} \uu\otimes \uu : \nabla^2 q \, \d x
+\int_{\T^2} \frac{\rho'(\phi)}{\rho(\phi)} ((\nabla \mu \cdot \nabla) \uu) \cdot \nabla q \, \d x- \int_{\T^2} \frac{\nu(\phi)}{\rho(\phi)} D \uu : \nabla^2 q \, \d x \\
&\quad + \int_{\T^2} \nu(\phi) D \uu: \Big( \frac{\rho'(\phi)}{\rho(\phi)^2}\nabla q \otimes \nabla \phi \Big) \, \d x- \int_{\T^2} \div(\nabla \phi\otimes \nabla \phi) \cdot \frac{\nabla q}{\rho(\phi)} \, \d x.
\end{split}
\end{equation}
Exploiting \eqref{LADY}, \eqref{Agmon2d}, \eqref{uL2H1-T} and \eqref{mu-H1eT}, we find
\begin{align*}
\begin{split}
\Big| \int_{\T^2} \uu\otimes \uu : \nabla^2 q \, \d x \Big|
&\leq \| \uu\|_{L^4(\T^2)}^2 \| q\|_{H^2(\Omega)}\\
&\leq C\|\uu \|_{H^1(\T^2)} (1+ \| \nabla \phi\|_{L^\infty(\T^2)}) \| P\|_{\LT},
\end{split}\\[10pt]
\begin{split}
\Big| \int_{\T^2} \frac{\rho'(\phi)}{\rho(\phi)} ((\nabla \mu \cdot \nabla) \uu) \cdot \nabla q \, \d x \Big| &\leq \Big\| \frac{\rho'(\phi)}{\rho(\phi)}\Big\|_{L^\infty(\T^2)}  \| \nabla \uu\|_{\LT} \| \nabla \mu \|_{L^\infty(\T^2)} \| \nabla q\|_{\LT}\\
&\leq  C \| \nabla \uu\|_{\LT} \| \nabla \mu \|_{\LT}^\frac12 \| \mu\|_{H^3(\T^2)}^\frac12 \| P\|_{\LT},
\end{split}\\[10pt]
\begin{split}
\Big|  - \int_{\T^2} \frac{\nu(\phi)}{\rho(\phi)} D \uu : \nabla^2 q \, \d x \Big| &\leq \Big\| \frac{\nu(\phi)}{\rho(\phi)}\Big\|_{L^\infty(\T^2)}  \| D\uu\|_{\LT} \| q\|_{H^2(\T^2)}\\
&\leq C \| D\uu\|_{\LT} (1+ \| \nabla \phi\|_{L^\infty(\T^2)}) \| P\|_{\LT},
\end{split}\\[10pt]
\begin{split}
\Big| \int_{\T^2} \nu(\phi) D \uu: \Big( \frac{\rho'(\phi)}{\rho(\phi)^2}\nabla q \otimes \nabla \phi \Big) \, \d x  \Big| &\leq \nu^\ast \| D\uu\|_{\LT} \Big\| \frac{\rho'(\phi)}{\rho(\phi)^2}\Big\|_{L^\infty(\T^2)} \| \nabla q\|_{\LT} \| \nabla \phi \|_{L^\infty(\T^2)}\\
&\leq C \| D\uu\|_{\LT} \| \nabla \phi\|_{L^\infty(\T^2)} \| P\|_{\LT},
\end{split}\\[10pt]
\begin{split}
\Big| - \int_{\T^2} \div(\nabla \phi\otimes \nabla \phi) \cdot \frac{\nabla q}{\rho(\phi)} \, \d x \Big| &\leq C \| \phi\|_{H^2(\T^2)} \| \nabla \phi\|_{L^\infty(\T^2)} \|\nabla q \|_{\LT}\\
&\leq C \| \phi\|_{H^2(\T^2)} \| \nabla \phi\|_{L^\infty(\T^2)} \| P\|_{\LT}.
\end{split}
\end{align*}
Thus, we are led to
\begin{equation}
\label{PL2}
\begin{split}
\| P\|_{\LT} 
&\leq C\|\uu \|_{H^1(\T^2)} (1+ \| \nabla \phi\|_{L^\infty(\T^2)})+ 
C \| \nabla \uu\|_{\LT}  \| \nabla \mu \|_{\LT}^\frac12 \| \mu\|_{H^3(\T^2)}^\frac12 \\
&\quad + C \| \phi\|_{H^2(\T^2)} \| \nabla \phi\|_{L^\infty(\T^2)}.
\end{split}
\end{equation}
Inserting \eqref{PL2} in \eqref{W4}, we obtain
\begin{equation}
\label{W4-2}
|W_4| \leq \sum_{i=1}^6 Y_i,
\end{equation}
where
\begin{align*}
&Y_1=C\|\uu \|_{H^1(\T^2)} (1+ \| \nabla \phi\|_{L^\infty(\T^2)})\| \phi\|_{H^2(\T^2)}^\frac12 (1+\| \nabla \mu\|_{L^2(\T^2)}^\frac34) \| \mu\|_{H^3(\T^2)}^\frac14, \\
&Y_2=C\|\uu \|_{H^1(\T^2)} (1+ \| \nabla \phi\|_{L^\infty(\T^2)}) (1+\| \nabla \mu\|_{\LT}^\frac12) \| \mu\|_{H^3(\T^2)}^\frac12,\\
&Y_3=C \| \nabla \uu\|_{\LT}  \| \phi\|_{H^2(\T^2)}^\frac12 (1+\| \nabla \mu\|_{L^2(\T^2)}^\frac54) \| \mu\|_{H^3(\T^2)}^\frac34,\\ 
&Y_4=C \| \nabla \uu\|_{\LT} (1+\| \nabla \mu\|_{\LT}) \| \mu\|_{H^3(\T^2)},\\
&Y_5=C \| \phi\|_{H^2(\T^2)}^\frac32 \| \nabla \phi\|_{L^\infty(\T^2)}(1+\| \nabla \mu\|_{L^2(\T^2)}^\frac34) \| \mu\|_{H^3(\T^2)}^\frac14,\\
&Y_6= C \| \phi\|_{H^2(\T^2)} \| \nabla \phi\|_{L^\infty(\T^2)}(1+\| \nabla \mu\|_{\LT}^\frac12) \| \mu\|_{H^3(\T^2)}^\frac12.
\end{align*}
By \eqref{phiH2eT}, \eqref{phiW2peT} and the Young inequality, we have
\begin{align}
\label{y1}
\begin{split}
&Y_1 
\begin{aligned}[t] 
&\leq \frac{\varepsilon}{48} \| \mu\|_{H^3(\T^2)}^2
+C \|\uu \|_{H^1(\T^2)}^\frac87 (1+ \| \nabla \phi\|_{L^\infty(\T^2)}^\frac87)\| \phi\|_{H^2(\T^2)}^\frac47 (1+\| \nabla \mu\|_{L^2(\T^2)}^\frac67)\\
&\leq \frac{\varepsilon}{48} \| \mu\|_{H^3(\T^2)}^2
+C \|\uu \|_{H^1(\T^2)}^\frac87 (1+ \| \nabla \phi\|_{L^\infty(\T^2)}^\frac87)(1+\| \nabla \mu\|_{L^2(\T^2)}^\frac87)\\
&\leq \frac{\varepsilon}{48} \| \mu\|_{H^3(\T^2)}^2
+C\big( 1+ \|\uu \|_{H^1(\T^2)}^2+ \| \nabla \phi\|_{L^\infty(\T^2)}^\frac83\big)(1+\| \nabla \mu\|_{L^2(\T^2)}^\frac87),
\end{aligned} 
\end{split}
\\[10pt]
\label{y2}
\begin{split}
&Y_2 
\begin{aligned}[t]
&\leq \frac{\varepsilon}{48} \| \mu\|_{H^3(\T^2)}^2
+C\|\uu \|_{H^1(\T^2)}^\frac43 (1+ \| \nabla \phi\|_{L^\infty(\T^2)}^\frac43) (1+\| \nabla \mu\|_{\LT}^\frac23)\\
&\leq \frac{\varepsilon}{48} \| \mu\|_{H^3(\T^2)}^2
+C\|\uu \|_{H^1(\T^2)}^\frac43 (1+ \| \phi\|_{W^{2,3}(\T^2)}^\frac43) (1+\| \nabla \mu\|_{\LT}^\frac23)\\
&\leq \frac{\varepsilon}{48} \| \mu\|_{H^3(\T^2)}^2
+C\|\uu \|_{H^1(\T^2)}^\frac43(1+\| \nabla \mu\|_{\LT}^2),
\end{aligned}
\end{split}
\\[10pt]
\label{y3}
\begin{split}
&Y_3
\begin{aligned}[t]
&\leq \frac{\varepsilon}{48} \| \mu\|_{H^3(\T^2)}^2
+ C \| \nabla \uu\|_{\LT}^\frac85  \| \phi\|_{H^2(\T^2)}^\frac45 (1+\| \nabla \mu\|_{L^2(\T^2)}^2) \\
&\leq  \frac{\varepsilon}{48} \| \mu\|_{H^3(\T^2)}^2
+ C \big( \| \nabla \uu\|_{\LT}^2 + \| \phi\|_{H^2(\T^2)}^4 \big) (1+\| \nabla \mu\|_{L^2(\T^2)}^2),
\end{aligned}
\end{split}
\\[10pt]
\label{y4}
\begin{split}
&Y_4
\begin{aligned}[t]
&\leq \frac{\varepsilon}{48} \| \mu\|_{H^3(\T^2)}^2
+ C \| \nabla \uu\|_{\LT}^2 (1+\| \nabla \mu\|_{L^2(\T^2)}^2),
\end{aligned}
\end{split}
\\[10pt]
\label{y5}
\begin{split}
&Y_5 
\begin{aligned}[t]
&\leq \frac{\varepsilon}{48} \| \mu\|_{H^3(\T^2)}^2
+
C \| \phi\|_{H^2(\T^2)}^\frac{12}{7} \| \nabla \phi\|_{L^\infty(\T^2)}^\frac87 
(1+\| \nabla \mu\|_{L^2(\T^2)}^\frac67)\\
&\leq \frac{\varepsilon}{48} \| \mu\|_{H^3(\T^2)}^2
+
C \| \nabla \phi\|_{L^\infty(\T^2)}^\frac{8}{7} 
(1+\| \nabla \mu\|_{L^2(\T^2)}^\frac{12}{7}),
\end{aligned}
\end{split}
\\[10pt]
\label{y6}
\begin{split}
&Y_6
\begin{aligned}[t]
&\leq 
\frac{\varepsilon}{48} \| \mu\|_{H^3(\T^2)}^2+
C \| \phi\|_{H^2(\T^2)}^\frac43 \| \nabla \phi\|_{L^\infty(\T^2)}^\frac43 
(1+\| \nabla \mu\|_{\LT}^\frac23) \\
&\leq 
\frac{\varepsilon}{48} \| \mu\|_{H^3(\T^2)}^2+
C  \| \nabla \phi\|_{L^\infty(\T^2)}^\frac43 
(1+\| \nabla \mu\|_{\LT}^\frac43).
\end{aligned}
\end{split}
\end{align}
Combining \eqref{W4-2} with \eqref{y1}-\eqref{y6}, we infer that 
\begin{equation}
\label{W4-3}
|W_4| \leq \frac{\varepsilon}{8} \| \mu\|_{H^3(\T^2)}^2+
C \big( 1+ \| \uu \|_{H^1(\T^2)}^2+\| \phi\|_{H^2(\T^2)}^4+
 \| \nabla \phi\|_{L^\infty(\T^2)}^\frac83 \big) 
(1+\| \nabla \mu\|_{\LT}^2).
\end{equation}
Collecting \eqref{W1}-\eqref{W3} and \eqref{W4-3} together, we find
\begin{align*}
\Big| \int_{\T^2}  \mu \partial_t \uu \cdot \nabla \phi \, \d x \Big|
\leq \frac{\varepsilon}{2} \| \mu\|_{H^3(\T^2)}^2+
C \big( 1+ \| \uu \|_{H^1(\T^2)}^2+\| \phi\|_{H^2(\T^2)}^4+
 \| \nabla \phi\|_{L^\infty(\T^2)}^\frac83 \big) 
(1+\| \nabla \mu\|_{\LT}^2).
\end{align*}
Hence, it follows from \eqref{HO3} and the above inequality that
\begin{equation}
\label{HO4}
\begin{split}
&\ddt \bigg[ \int_{\T^2} \frac12 |\nabla \mu|^2 \, \d x 
+ \int_{\T^2} \mu (\uu \cdot \nabla \phi)\, \d x \bigg] 
+ \frac14 \| \nabla \partial_t \phi\|_{\LT}^2+ \frac{\varepsilon}{2} \| \mu\|_{H^3(\T^2)}^2\\
&\quad 
\leq C \big( 1+ \| \uu \|_{H^1(\T^2)}^2+\| \phi\|_{H^2(\T^2)}^4+
 \| \nabla \phi\|_{L^\infty(\T^2)}^\frac83 \big) 
(1+\| \nabla \mu\|_{\LT}^2).
\end{split}
\end{equation}
We now set 
$$
X(t)= \int_{\T^2} \frac12 |\nabla \mu(t)|^2 \, \d x 
+ \int_{\T^2} \mu(t) (\uu(t) \cdot \nabla \phi(t)) \, \d x .
$$
Thanks to \eqref{uL2H1-T}, we observe that 
$$
\int_{\T^2} \mu (\uu  \cdot \nabla \phi)
= \int_{\T^2} \phi( \uu\cdot \nabla \mu) \, \d x\\
\leq \| \phi\|_{L^\infty(\T^2)} \| \uu\|_{\LT} \|\nabla \mu\|_{\LT}
\leq C \| \nabla \mu\|_{\LT}.
$$
Then, there exists a positive constant $\overline{C}$ depending on $E(\uu_0,\phi_0)$ such that
$$
X(t)\geq \frac14 \| \nabla \mu(t)\|_{\LT}^2 -\overline{C}.
$$
Therefore, we deduce the differential inequality
\begin{equation}
\label{HO42}
\ddt X(t)
+ \frac14 \| \nabla \partial_t \phi\|_{\LT}^2+ \frac{\varepsilon}{2} \| \mu\|_{H^3(\T^2)}^2
\leq Y(t) (1+\overline{C}+ X(t)),
\end{equation}
where
$$
Y(t)= C\big(1+ \| \uu(t) \|_{H^1(\T^2)}^2+\| \phi(t)\|_{H^2(\T^2)}^4+
 \| \nabla \phi (t)\|_{L^\infty(\T^2)}^\frac83 \big).
$$
In light of \eqref{RSS-T} and \eqref{phi-T}, we infer from the Gagliardo-Nirenberg inequality \eqref{GN-Linf} with $s=3$ that $\| \phi \|_{L^\frac83(0,T;W^{1,\infty}(\T^2))}\leq C(1+T)$, for all $T<T_\ast$.
In turn, it gives $\| Y\|_{L^1(0,T)}\leq C(1+T)$ for all $T<T_\ast$ (cf. \eqref{uL2H1-T} and \eqref{phi-T}). Thus, the Gronwall lemma yields
\begin{equation}
\sup_{t\in [0,T]} X(t) \leq \Big(X(0)+ (1+\overline{C}) \int_0^{T} Y(\tau)\, \d \tau \Big) \mathrm{e}^{\int_0^{T} Y(\tau)\, \d \tau}, \quad \forall \, T<T_\ast,
\end{equation}
which entails that
\begin{equation}
\label{muH1-GT}
\sup_{t\in [0,T]} \| \mu(t)\|_{H^1(\T^2)} \leq K_T, \quad \forall \, T<T_\ast,
\end{equation}
where $K_T$ stands for a generic constant depending on the parameters of the system, the initial energy $E(\uu_0,\phi_0)$, the norms of the initial data $\| \uu_0\|_{H^1(\T^2)}$ and $\| \mu_0\|_{H^1(\T^2)}$, and the time $T$. In particular, $K_T$ is finite for any $T<\infty$.
Integrating in time \eqref{HO42}, we infer that
\begin{equation}
\label{muH3-GT}
\int_0^{T} \| \partial_t \phi (\tau)\|_{H^1(\T^2)}^2
+ \| \mu (\tau)\|_{H^3(\T^2)}^2 \, \d \tau \leq K_T, \quad \forall \, T<T_\ast.
\end{equation}
As a consequence, we obtain from \eqref{phiW2peT} that
\begin{equation}
\label{phiW2p-GT}
\sup_{t\in [0,T]} \| \phi(t)\|_{W^{2,p}(\T^2)} 
+ \sup_{t\in [0,T]} \| F'(\phi(t))\|_{L^p(\T^2)} \leq K_{p,T}, \quad \forall \, T<T_\ast, \ \forall \, p \in [2,\infty).
\end{equation}
Finally, as in Section \ref{Loc-Ex},  by exploiting \cite[Lemma A.6]{CG2020}, we immediately deduce that
\begin{equation}
\sup_{t\in[0,T]} \| F''(\phi(t))\|_{L^p(\T^2)}
+ \sup_{t\in[0,T]} \| F'''(\phi(t))\|_{L^p(\T^2)}
\leq K_{p,T}, \quad \forall \, T<T_\ast, \ \forall \, p \in [2,\infty),
\end{equation}
which implies that 
\begin{equation*}
\sup_{t\in [0,T]} \| \phi(t)\|_{H^3(\T^2)}+ \int_0^{T} \| \partial_t \mu (\tau)\|_{(H^1(\T^2))'}^2 \, \d \tau\leq K_{T}, \quad \forall \, T<T_\ast.
\end{equation*}

\subsection{High-Order Estimates for the Velocity Field}
Multiplying \eqref{AGG-T} by $\partial_t \uu$ and integrating over $\T^2$, we obtain (cf. \eqref{High1})
\begin{equation}
\label{HO5}
\begin{split}
&\frac12 \ddt \int_{\T^2} \nu(\phi) |D \uu|^2 \, \d x
+ \int_{\T^2} \rho(\phi) |\partial_t \uu|^2 \, \d x\\
&\quad =- \int_{\T^2} \rho(\phi) ( (\uu \cdot \nabla) \uu ) \cdot \partial_t \uu \, \d x + \frac{\nu_1-\nu_2}{2}  \int_{\T^2}\partial_t \phi  |D \uu|^2 \, \d x\\
&\qquad +  \frac{\rho_1-\rho_2}{2} \int_{\T^2} ((\nabla \mu\cdot \nabla) \uu) \cdot \partial_t \uu - \int_{\T^2} \Delta \phi \nabla \phi \cdot \partial_t \uu \, \d x.
\end{split}
\end{equation}
On the other hand, multiplying \eqref{AGG-T} by $-\Delta \uu$, we find
\begin{equation}
\label{HO6}
\begin{split}
\frac{\nu_\ast}{2} \| \Delta \uu\|_{\LT}^2
&\leq  \int_{\T^2} \rho(\phi) \partial_t \uu \cdot \Delta \uu \, \d x
+ \int_{\T^2} \rho(\phi) ((\uu \cdot \nabla)\uu) \cdot \Delta \uu \, \d x\\
&\quad -\int_{\T^2} \rho'(\phi) ((\nabla \mu \cdot \nabla) \uu) \cdot \Delta \uu \, \d x
+\int_{\T^2} \Delta \phi \nabla \phi \cdot \Delta \uu \, \d x\\
&\quad - \int_{\T^2} \nu'(\phi) (D\uu \nabla \phi) \cdot \Delta \uu \, \d x.
\end{split}
\end{equation}
By the Young inequality, we simply have
\begin{equation}
\label{HO7}
\begin{split}
\frac{\nu_\ast}{4} \| \Delta \uu\|_{\LT}^2
&\leq  \frac{(\rho^\ast)^2}{\nu_\ast} \| \partial_t \uu \|_{\LT}^2
+ \int_{\T^2} \rho(\phi) ((\uu \cdot \nabla)\uu) \cdot \Delta \uu \, \d x\\
&\quad -\int_{\T^2} \rho'(\phi)((\nabla \mu\cdot \nabla) \uu) \cdot \Delta \uu \, \d x
+\int_{\T^2} \Delta \phi \nabla \phi \cdot \Delta \uu \, \d x\\
&\quad - \int_{\T^2} \nu'(\phi) (D\uu \nabla \phi) \cdot \Delta \uu \, \d x.
\end{split}
\end{equation}
Multiplying \eqref{HO7} by $\frac{\nu_\ast}{2\rho^\ast}$ and adding the resulting inequality to \eqref{HO5}, we reach
\begin{equation}
\label{HO8}
\begin{split}
&\frac12 \ddt \int_{\T^2} \nu(\phi) |D \uu|^2 \, \d x
+ \frac{\rho_\ast}{2} \|\partial_t \uu\|_{\LT}^2 + \frac{\nu_\ast^2}{8\rho^\ast}
\| \Delta \uu\|_{\LT}^2\\
&\quad \leq - \int_{\T^2} \rho(\phi) \big( (\uu \cdot \nabla) \uu \big) \cdot \partial_t \uu \, \d x +  \frac{\nu_1-\nu_2}{2} \int_{\T^2} \partial_t \phi  |D \uu|^2 \, \d x\\
&\qquad +  \frac{\rho_1-\rho_2}{2} \int_{\T^2} ((\nabla \mu \cdot \nabla) \uu) \cdot \partial_t \uu \, \d x- \int_{\T^2} \Delta \phi \nabla \phi \cdot \partial_t \uu \, \d x \\
&\qquad+\frac{\nu_\ast}{2\rho^\ast}  \int_{\T^2} \rho(\phi) ((\uu \cdot \nabla)\uu) \cdot \Delta \uu \, \d x-\frac{\nu_\ast}{2\rho^\ast} \int_{\T^2} \rho'(\phi)((\nabla \mu \cdot \nabla) \uu) \cdot \Delta \uu \, \d x \\
&\qquad 
+ \frac{\nu_\ast}{2\rho^\ast} \int_{\T^2} \Delta \phi \nabla \phi \cdot \Delta \uu \, \d x -\frac{\nu_\ast}{2\rho^\ast}  \int_{\T^2} \nu'(\phi) (D\uu \nabla \phi) \cdot \Delta \uu \, \d x\\
&\quad =\sum_{i}^8 L_i.
\end{split}
\end{equation}
Notice that 
$$
\| \uu\|_{H^1(\T^2)}\leq C (1+ \| D \uu\|_{\LT}), \quad
\| \uu\|_{H^2(\T^2)}\leq C(1+\| \Delta \uu\|_{\LT})
$$ 
due to \eqref{KORN2}, \eqref{H2equiv-T} and \eqref{uL2H1-T}. 
By \eqref{LADY}, \eqref{uL2H1-T}, \eqref{phiW2p-GT}, we can estimate the terms $L_i$ as follows 
\begin{align}
\label{L1}
\begin{split}
&L_1
\begin{aligned}[t]
&\leq C \| \uu\|_{L^4(\T^2)} \| \nabla \uu\|_{L^4(\T^2)} \| \partial_t \uu\|_{\LT}\\
&\leq C \| \uu \|_{H^1(\T^2)} (1+ \| \Delta \uu\|_{\LT}^\frac12)
\| \partial_t \uu\|_{\LT}\\
&\leq \frac{\rho_\ast}{12} \| \partial_t \uu\|_{\LT}^2 +
\frac{\nu_\ast^2}{96\rho^\ast} \| \Delta \uu\|_{\LT}^2 + C\| D \uu\|_{\LT}^2(1+ \| D \uu\|_{\LT}^2),
\end{aligned}
\end{split}
\\[10pt]
\label{L2}
\begin{split}
&L_2
\begin{aligned}[t]
&\leq C \| \partial_t \phi \|_{L^6(\T^2)} \| D \uu \|_{\LT}  \|D \uu \|_{L^3(\T^2)}\\
&\leq C\| \partial_t \phi \|_{H^1(\T^2)} \| D \uu \|_{\LT}  
(1+ \|\Delta \uu \|_{\LT})\\
&\leq \frac{\nu_\ast^2}{96\rho^\ast} \| \Delta \uu\|_{\LT}^2 +  C(1 +\| \partial_t \phi \|_{H^1(\T^2)}^2)(1+ \| D \uu \|_{\LT}^2),
\end{aligned}
\end{split}
\\[10pt]
\label{L3}
\begin{split}
&L_3
\begin{aligned}[t]
&\leq 
C \| \nabla \uu\|_{\LT} \| \nabla \mu\|_{L^\infty(\T^2)} 
\| \partial_t \uu\|_{\LT}\\
&\leq \frac{\rho_\ast}{12} \| \partial_t \uu\|_{\LT}^2
+C\| \mu\|_{H^3(\T^2)}^2 \| D \uu\|_{\LT}^2,
\end{aligned}
\end{split}
\\[10pt]
\label{L4}
\begin{split}
&L_4
\begin{aligned}[t]
&\leq \| \Delta \phi\|_{\LT} \| \nabla \phi\|_{L^\infty(\T^2)} \| \partial_t \uu\|_{\LT}\\
&\leq \frac{\rho_\ast}{12} \| \partial_t \uu\|_{\LT}^2
+K_T,
\end{aligned}
\end{split}
\\[10pt]
\label{L5}
\begin{split}
&L_5 
\begin{aligned}[t]
&\leq C \| \uu \|_{L^4(\T^2)} \| \nabla \uu\|_{L^4(\T^2)} \| \Delta\uu\|_{\LT}\\
&\leq C (1+ \| D \uu\|_{\LT}^\frac12)\| D \uu\|_{\LT}^\frac12 (1+ \| \Delta \uu\|_{\LT}^\frac12)\| \Delta\uu \|_{\LT} \\
&\leq \frac{\nu_\ast^2}{96\rho^\ast} \| \Delta \uu\|_{\LT}^2+C (1+\| D \uu\|_{\LT}^2)^2,
\end{aligned}
\end{split}
\\[10pt]
\label{L6}
\begin{split}
&L_6
\begin{aligned}[t]
&\leq 
C \| \nabla \uu\|_{\LT} \| \nabla \mu\|_{L^\infty(\T^2)}\| \Delta \uu\|_{\LT}\\
&\leq \frac{\nu_\ast^2}{96\rho^\ast} \| \Delta \uu\|_{\LT}^2+ C \| \mu\|_{H^3(\T^2)}^2 \| D \uu\|_{\LT}^2,
\end{aligned}
\end{split}
\\[10pt]
\label{L7}
\begin{split}
&L_7
\begin{aligned}[t]
&\leq C \|\Delta \phi \|_{\LT}\| \nabla \phi\|_{L^\infty(\T^2)} \| \Delta \uu\|_{\LT}\\
&\leq \frac{\nu_\ast^2}{96\rho^\ast} \| \Delta \uu\|_{\LT}^2+ K_T,
\end{aligned}
\end{split}
\\[10pt]
\label{L8}
\begin{split}
&L_8
\begin{aligned}[t]
&\leq C \| D \uu\|_{\LT} \| \nabla \phi\|_{L^\infty(\T^2)} \| \Delta \uu\|_{\LT}\\
&\leq \frac{\nu_\ast^2}{96\rho^\ast} \| \Delta \uu\|_{\LT}^2+ K_T \| D \uu\|_{\LT}^2. 
\end{aligned}
\end{split}
\end{align}
Hence, it follows that on $[0,T]$, for all $T<T_\ast$,
\begin{equation*}
\frac12 \ddt \int_{\T^2} \nu(\phi) |D \uu|^2 \, \d x
+ \frac{\rho_\ast}{4} \|\partial_t \uu\|_{\LT}^2 
+ \frac{\nu_\ast^2}{16 \rho^\ast}
\| \Delta \uu\|_{\LT}^2\leq \widetilde{Y}(t)
(1+ \| D \uu\|_{\LT}^2),
\end{equation*}
where
$$
\widetilde{Y}(t)= K_T \big(1+ \| D \uu(t)\|_{\LT}^2+ \| \partial_t \phi(t)\|_{H^1(\T^2)}^2+
\|\mu(t) \|_{H^3(\T^2)}^2 \big).
$$
In light of \eqref{uL2H1-T} and \eqref{muH3-GT}, an application of the Gronwall lemma yields
\begin{equation}
\label{uH1-GT}
\sup_{t\in[0,T]} \| \uu(t)\|_{H^1(\T^2)}^2
 + \int_0^{T}\|\partial_t \uu(\tau)\|_{\LT}^2 
+ \| \Delta \uu (\tau)\|_{\LT}^2 \, \d \tau
\leq \widetilde{K}_T, \quad \forall \, T <T_\ast,
\end{equation}
where 
$$
\widetilde{K}_T = C \bigg( \| \uu_0\|_{H^1(\T^2)}^2+ \Big( \int_0^T \widetilde{Y}(\tau)\, \d \tau\Big)^2 \bigg) \mathrm{e}^{\int_0^T \widetilde{Y}(\tau)\, \d \tau}
$$
for some positive constant $C$ depending on $\nu_\ast$, $\rho_\ast$ and $\rho^\ast$.

\subsection{Global Existence of Strong Solutions} 
The uniform-in-time estimates \eqref{muH1-GT}-\eqref{phiW2p-GT} and \eqref{uH1-GT} entails that the solution does not blowup as $T$ approaches $T_\ast$. 
By a classical argument, it is possible to show that $(\uu(T),\mu(T))$ is a Cauchy sequence in $\VV_\sigma \times H^1(\T^2)$ as $T \rightarrow T_\ast$. This implies that $\uu(T_\ast)$ and $\mu(T_\ast)$ are well-defined. Then, the solution can be continued beyond $T_\ast$ into a solution which satisfies \eqref{RSS-T} and \eqref{AGG-T} on an interval $(0,\overline{T})$ for some $\overline{T}>T_\ast$. This contradicts the maximality of $T_\ast$. Hence, $T_\ast =\infty$.

\section{Uniqueness}
\label{Un}
\setcounter{equation}{0}
In this section we show the continuous dependence on the initial data and the uniqueness of the strong solutions proved in Theorem \ref{mr1} and Theorem \ref{mr2}. We demonstrate hereafter the case of a general bounded domain $\Omega\subset \mathbb{R}^2$. The proof in the case $\Omega=\mathbb{T}^2$ can be adapted with minor changes.
\smallskip

Let $(\uu_1,P_1,\phi_1)$ and $(\uu_2,P_2,\phi_2)$ be two strong solutions to system \eqref{AGG} with boundary conditions \eqref{AGG-bc} defined on a common interval $[0,T_0]$ given by Theorem \ref{mr1}. We consider $\uu=\uu_1-\uu_2$, $P=P_1-P_2$ and $\phi=\phi_1-\phi_2$. It is clear that
\begin{align}
\label{U-u}
\begin{split}
&\rho(\phi_1)\partial_t \uu + (\rho(\phi_1)-\rho(\phi_2)) \partial_t \uu_2 +
\big(\rho(\phi_1)(\uu_1 \cdot \nabla) \uu_1- \rho(\phi_2)(\uu_2 \cdot \nabla) \uu_2\big)\\
&\quad - \frac{\rho_1-\rho_2}{2}\big( (\nabla \mu_1\cdot \nabla) \uu_1-(\nabla \mu_2 \cdot \nabla) \uu_2 \big)
- \div (\nu(\phi_1)D\uu) 
-\div( (\nu(\phi_1)-\nu(\phi_2))D \uu_2)\\
&\quad + \nabla P= - \div(\nabla \phi_1 \otimes \nabla \phi_1+ \nabla \phi_2\otimes \nabla \phi_2),
\end{split}
\\[10pt]
\label{U-phi}
\begin{split}
&\partial_t \phi +\uu_1\cdot \nabla \phi + \uu \cdot \nabla \phi_2= \Delta \mu,\\
&\mu= -\Delta \phi+\Psi'(\phi_1)- \Psi'(\phi_2),
\end{split}
\end{align}
almost everywhere in $\Omega \times (0,T_0)$.
Multiplying \eqref{U-u} by $\uu$ and integrating over $\Omega$, we find
\begin{equation}
\label{U-1}
\begin{split}
&\frac12 \ddt \int_{\Omega} \rho(\phi_1) |\uu|^2 \, \d x 
+ \int_{\Omega} \nu(\phi_1) |D \uu|^2 \, \d x
\\
&= - \int_{\Omega} (\rho(\phi_1)-\rho(\phi_2) ) \partial_t \uu_2 \cdot \uu \, \d x
-\int_{\Omega} \rho(\phi_1) (\uu \cdot \nabla) \uu_2 \cdot \uu \, \d x\\
&\quad 
- \int_{\Omega} (\rho(\phi_1)-\rho(\phi_2))(\uu_2\cdot \nabla)\uu_2 \cdot \uu\, \d x
+\frac{\rho_1-\rho_2}{2} \int_{\Omega} ((\nabla \mu \cdot\nabla) \uu_2)  \cdot \uu \, \d x \\
&\quad -\int_{\Omega} (\nu(\phi_1)-\nu(\phi_2))D \uu_2 : \nabla \uu \, \d x
+ \int_{\Omega} (\nabla \phi_1 \otimes \nabla \phi + \nabla \phi\otimes \nabla \phi_2): \nabla \uu \, \d x\\
&= \sum_{i=1}^6 Z_i.
\end{split}
\end{equation} 
Here we have used that 
$$
-\int_{\Omega} \partial_t \rho(\phi_1) \frac{|\uu|^2}{2}\, \d x 
+\int_{\Omega} \rho(\phi_1) (\uu_1\cdot \nabla)\uu \cdot \uu \, \d x 
-\frac{\rho_1-\rho_2}{2} \int_{\Omega} (\nabla \mu_1\cdot \nabla) \uu \cdot \uu \, \d x=0.  
$$
Taking the gradient of \eqref{U-phi}, multiplying the resulting equation by $\nabla \phi$ and integrating over $\Omega$, then using the boundary conditions \eqref{AGG-bc}, we obtain
\begin{equation}
\label{U-2}
\begin{split}
&\frac12 \ddt \|\nabla \phi\|_{\L2}^2 +  \| \nabla \Delta \phi\|_{\L2}^2\\
&\quad = - \int_{\Omega} \nabla (\uu_1\cdot \nabla \phi) \cdot \nabla \phi \, \d x -\int_{\Omega} \nabla(\uu \cdot \nabla \phi_2) \cdot \nabla \phi \, \d x\\
&\qquad+\int_{\Omega} \nabla(\Psi'(\phi_1)-\Psi'(\phi_2)) \cdot \nabla \Delta \phi \, \d x\\
&\quad = \sum_{i=7}^{9} Z_i.
\end{split}
\end{equation}
Since 
$$
\ddt \overline{\phi}=0,
$$
by \eqref{U-1} and \eqref{U-2}, we reach
$$
\frac12 \ddt \bigg[ \int_{\Omega} \rho(\phi_1) |\uu|^2 \, \d x +  \|\nabla \phi\|_{\L2}^2 + |\overline{\phi}|^2 \bigg]
+ \int_{\Omega} \nu(\phi_1) |D \uu|^2 \, \d x
+  \| \nabla \Delta \phi\|_{\L2}^2= 
\sum_{i=1}^{9} Z_i.
$$
We recall that $\| \phi\|_{H^3(\Omega)} \leq C\big(\| \phi\|_{H^1(\Omega)} + \| \nabla \Delta \phi\|_{\L2})$. 
By exploiting \eqref{LADY}, \eqref{KORN} and the regularity of the strong solutions, we infer that
\begin{align}
\label{Z1}
\begin{split}
&|Z_1|
\begin{aligned}[t]
&\leq C\| \phi\|_{L^6(\Omega)}\| \partial_t \uu_2\|_{\L2} \| \uu\|_{L^3(\Omega)}\\
&\leq \frac{\nu_\ast}{8}\| D \uu\|_{\L2}^2
+C \| \partial_t \uu_2\|_{\L2}^2 \big(\| \uu \|_{\L2}^2 + \| \phi \|_{H^1(\Omega)}^2  \big),
\end{aligned}
\end{split}
\\[10pt]
\label{Z2}
\begin{split}
&|Z_2|
\begin{aligned}[t]
&\leq C \| \uu\|_{L^3(\Omega)}\| \nabla \uu_2\|_{L^6(\Omega)} \| \uu \|_{\L2}\\
&\leq  \frac{\nu_\ast}{8}\| D \uu\|_{\L2}^2+C \| \uu_2\|_{H^2(\Omega)}^2
\|\uu \|_{\L2}^2,
\end{aligned}
\end{split}
\\[10pt]
\label{Z3}
\begin{split}
&|Z_3|
\begin{aligned}[t]
&\leq C \| \phi\|_{L^6(\Omega)} \| \uu_2\|_{L^6(\Omega)}  \|\nabla \uu_2\|_{L^6(\Omega)}  \| \uu \|_{\L2}\\
&\leq C\|\nabla \uu_2\|_{L^6(\Omega)} \big(  \| \uu \|_{\L2}^2 + \| \phi \|_{H^1(\Omega)}^2 \big),
\end{aligned}
\end{split}
\\[10pt]
\label{Z4}
\begin{split}
&|Z_4| 
\begin{aligned}[t]
&\leq \| \nabla \uu_2\|_{L^4(\Omega)}   \| \nabla \mu\|_{L^2} \| \uu\|_{L^4(\Omega)}\\
&\leq C \| \uu_2\|_{H^2(\Omega)}^\frac12 \big( \| \nabla \Delta \phi\|_{\L2}+
\| \Psi''(\phi_1) \nabla \phi\|_{\L2}+ \| (\Psi''(\phi_1)-\Psi''(\phi_2))\nabla \phi_2\|_{\L2} \big) \\
&\quad \times  \| \uu\|_{\L2}^\frac12 \| D \uu\|_{\L2}^\frac12\\
&\leq C \| \uu_2\|_{H^2(\Omega)}^\frac12 \big( \| \nabla \Delta \phi\|_{\L2}+
\| \Psi''(\phi_1) \|_{\L2} \| \nabla \phi\|_{L^\infty(\Omega)}\\
&\quad + (\| \Psi'''(\phi_1)\|_{\L2}+\|\Psi'''(\phi_2)\|_{\L2}) \| \phi\|_{L^\infty(\Omega)} \|\nabla \phi_2\|_{L^\infty(\Omega)} \big) \| \uu\|_{\L2}^\frac12 \| D \uu\|_{\L2}^\frac12\\
&\leq C \| \uu_2\|_{H^2(\Omega)}^\frac12 \big(\| \phi\|_{H^1(\Omega)} + \| \nabla \Delta \phi\|_{\L2} \big) \| \uu\|_{\L2}^\frac12 \| D \uu\|_{\L2}^\frac12\\
&\leq \frac16 \| \nabla \Delta \phi\|_{\L2}^2+
\frac{\nu_\ast}{8} \| D \uu\|_{\L2}^2+
C\| \uu_2\|_{H^2(\Omega)}^2 \big( \| \uu \|_{\L2}^2 + \| \phi \|_{H^1(\Omega)}^2 \big),
\end{aligned}
\end{split}
\\[10pt]
\label{Z5}
\begin{split}
&|Z_5| 
\begin{aligned}[t]
&\leq C \| \phi\|_{L^3(\Omega)} \| D \uu_2\|_{L^6(\Omega)} \| \nabla \uu\|_{\L2}\\
&\leq \frac{\nu_\ast}{8} \| D \uu\|_{\L2}^2+ C \| D \uu_2\|_{L^6(\Omega)}^2 \| \phi\|_{H^1(\Omega)}^2,
\end{aligned}
\end{split}
\\[10pt]
\label{Z6}
\begin{split}
&|Z_6|
\begin{aligned}[t]
&\leq C (\| \nabla \phi_1\|_{L^\infty(\Omega)}+\| \nabla \phi_2\|_{L^\infty(\Omega)}) \| \nabla \phi\|_{\L2} \| \nabla \uu\|_{\L2}\\
&\leq \frac{\nu_\ast}{8} \| D \uu\|_{\L2}^2+ C \| \phi\|_{H^1(\Omega)}^2,
\end{aligned}
\end{split}
\\[10pt]
\label{Z7}
\begin{split}
&|Z_7|
\begin{aligned}[t]
&\leq \big( \| \nabla \uu_1\|_{L^2(\Omega)}  \| \nabla \phi \|_{L^\infty(\Omega)} +
\| \uu_1\|_{L^6(\Omega)} \| \phi\|_{W^{2,3}(\Omega)} \big) \| \nabla \phi\|_{\L2}\\
&\leq C \big( \| \phi \|_{H^1(\Omega)}+ \| \nabla \Delta \phi\|_{\L2} \big) \|\nabla \phi \|_{\L2}\\
&\leq \frac16 \| \nabla \Delta \phi\|_{\L2}^2+
C\| \phi \|_{H^1(\Omega)}^2,
\end{aligned}
\end{split}
\\[10pt]
\label{Z8}
\begin{split}
&|Z_8| 
\begin{aligned}[t]
&\leq \big( \| \nabla \uu\|_{\L2} \| \nabla \phi_2 \|_{L^\infty(\Omega)}
+ \| \uu\|_{L^4(\Omega)} \| \phi_2 \|_{W^{2,4}(\Omega)} \big) \|\nabla \phi \|_{\L2}\\
&\leq \frac{\nu_\ast}{8} \| D \uu\|_{\L2}^2+C\| \nabla \phi\|_{\L2}^2,
\end{aligned}
\end{split}
\\[10pt]
\label{Z9}
\begin{split}
&|Z_9| 
\begin{aligned}[t]
&\leq 
 \| \Psi''(\phi_1) \nabla \phi\|_{\L2}\| \nabla \Delta \phi\|_{\L2}
 + \| (\Psi''(\phi_1)-\Psi''(\phi_2))\nabla \phi_2\|_{\L2} \| \nabla \Delta \phi\|_{\L2} \\
&\leq C \| \Psi''(\phi_1) \|_{L^6(\Omega)} \| \nabla \phi\|_{L^3(\Omega)}
\| \nabla \Delta \phi\|_{\L2}\\
&\quad + \big(\| \Psi'''(\phi_1)\|_{L^6(\Omega)}+\|\Psi'''(\phi_2)\|_{L^6(\Omega)} \big) 
\| \phi\|_{L^3(\Omega)} \|\nabla \phi_2\|_{L^\infty(\Omega)} \| \nabla \Delta \phi\|_{\L2}\\
&\leq C \| \phi\|_{H^2(\Omega)} \| \nabla \Delta \phi\|_{\L2}\\
&\leq C \| \phi\|_{H^1(\Omega)}^\frac12 \big( \| \phi\|_{H^1(\Omega)}+ \| \nabla \Delta \phi\|_{\L2}\big)^\frac12 
\| \nabla \Delta \phi\|_{\L2}\\
&\leq \frac16 \| \nabla \Delta \phi\|_{\L2}^2+ C \| \phi\|_{H^1(\Omega)}^2.
\end{aligned}
\end{split}
\end{align}
Therefore, by \eqref{Z1}-\eqref{Z9}, we find the differential inequality
\begin{equation*}
\begin{split}
&\frac12 \ddt \bigg[ \int_{\Omega} \rho(\phi_1) |\uu|^2 \, \d x +  \|\nabla \phi\|_{\L2}^2 + |\overline{\phi}|^2 \bigg]
+ \frac{\nu_\ast}{4}\|D \uu\|_{\L2}^2 
+ \frac12  \| \nabla \Delta \phi\|_{\L2}^2\\
&\qquad \leq \overline{C} \big( 1+ \| \partial_t \uu_2\|_{\L2}^2+
\| \uu_2\|_{H^2(\Omega)}^2  \big)  \big( \| \uu\|_{\L2}^2 + \|  \phi\|_{H^1(\Omega)}^2 \big),
\end{split}
\end{equation*}
where the constant $\overline{C}$ depends on the norm of the initial data and $T_0$.
Thanks to the Gronwall lemma, together with \eqref{normH1-2}, we deduce for all $t \in [0,T_0]$ that
$$
\| \uu (t)\|_{\L2}^2+ \| \phi(t)\|_{H^1(\Omega)}^2 
\leq C \big( \| \uu (0)\|_{\L2}^2+ \| \phi (0) \|_{H^1(\Omega)}^2 \big) \mathrm{e}^{\overline{C}\int_0^{T_0} (1+ \| \partial_t \uu_2(\tau)\|_{\L2}^2+
\| \uu_2 (\tau)\|_{H^2(\Omega)}^2) \, \d \tau}.
$$
The above inequality proves the continuous dependence of the solutions on the initial data. In particular, when $\uu(0)=\mathbf{0}$ and $\phi(0)=0$, it follows that $\uu(t)=\mathbf{0}$ and $\phi(t)=0$ for all $t\in [0,T_0]$. Thus, the strong solution is unique.


\begin{thebibliography}{99} \itemsep=2pt

\bibitem{A2009}
{\au H. Abels},
{\ti On a diffuse interface model for two-phase flows of viscous, incompressible fluids with matched densities},
{\jou Arch.\ Ration.\ Mech.\ Anal.}
\no{194}{463--506}{2009}

\bibitem{A2009-2}
{\au H. Abels},
{\ti Existence of weak solutions for a diffuse interface model for viscous, incompressible fluids with general densities},
{\jou Comm.\ Math.\ Phys.}
\no{289}{45--73}{2009}

\bibitem{A2012}
{\au H. Abels},
{\ti Strong well-posedness of a diffuse interface model for a viscous, quasi-incompressible two-phase flow},
{\jou SIAM J.\ Math.\ Anal.}
\no{44}{316--340}{2012}

\bibitem{AB2018}
{\au H. Abels, D. Breit}, 
{\ti Weak solutions for a non-Newtonian diffuse interface model with different densities},
{\jou Nonlinearity}
\no{29}{3426--3453}{2016}

\bibitem{ADG2013}
{\au H. Abels, D. Depner, H. Garcke},
{\ti Existence of weak solutions for a diffuse interface model for two-phase flows of incompressible fluids with different densities},
{\jou J.\ Math.\ Fluid Mech.}
\no{15}{453--480}{2013}

\bibitem{ADG2013-2}
{\au H. Abels, D. Depner, H. Garcke},
{\ti On an incompressible Navier-Stokes/Cahn-Hilliard system with degenerate mobility},
{\jou Ann.\ Inst.\ H.\ Poincar\'{e} Anal. Non Lin\'{e}aire}
\no{30}{1175--1190}{2013}

\bibitem{AGG2012}
{\au H. Abels, H. Garcke, G. Gr\"{u}n},
{\ti Thermodynamically consistent, frame indifferent diffuse interface models for incompressible two-phase flows with different densities},
{\jou Math.\ Models Methods Appl.\ Sci.}
\no{22}{1150013}{2012}

\bibitem{AG2018}
{\au H. Abels, H. Garcke}, 
{\bk Weak solutions and diffuse interface models for incompressible two-phase flows}.
\eds{Handbook of Mathematical Analysis in Mechanics of Viscous Fluid}{Springer International Publishing}{2018}

\bibitem{AT2020}
{\au H. Abels, Y. Terasawa},
{\ti Weak solutions for a diffuse interface model for two-phase flows of incompressible fluids with different densities and nonlocal free energies},
{\jou Math.\ Meth.\ Appl.\ Sci.}
\no{43}{3200--3219}{2020}

\bibitem{B2001}
{\au F. Boyer},
{\ti Nonhomogeneous Cahn-Hilliard fluids},
{\jou Ann.\ Inst.\ H.\ Poincar\'{e} Anal.\ Non Lin\'{e}aire}
\no{18}{225--259}{2001}

\bibitem{B2002}
{\au F. Boyer}, 
{\ti A theoretical and numerical model for the study of incompressible mixture flows}, 
{\jou Comput.\ Fluids} 
\no{31}{41--68}{2002}

\bibitem{CK2003}
{\au H. J. Choe, H. Kim}, 
{\ti Strong solutions of the Navier-Stokes equations for nonhomogeneous
incompressible fluids},
{\jou Comm.\ Partial Differential Equations} 
\no{28}{1183--1202}{2003}

\bibitem{CG2020}
{\au M. Conti, A. Giorgini},
{\ti Well-posedness for the Brinkman-Cahn-HIlliard system with unmatched viscosities},
{\jou J.\ Differential Equations}
\no{268}{6350--6384}{2020}


\bibitem{DSS2007}
{\au H. Ding, P.D.M. Spelt, C. Shu}, 
{\ti Diffuse interface model for incompressible two-phase
flows with large density ratios}, 
{\jou J.\ Comput.\ Phys.} 
\no{226}{2078--2095}{2007}

\bibitem{Frigeri2016}
{\au S. Frigeri},
{\ti Global existence of weak solutions for a nonlocal model for two-phase flows of incompressible fluids with unmatched densities}, {\jou Math.\ Models Methods Appl. Sci.} 
\no{26}{1957--1993}{2016}

\bibitem{Galdi}
{\au G.P. Galdi},
{\bk An Introduction to the Mathematical Theory of the Navier–Stokes Equations},
\eds{Vol. 1. Springer}{Berlin}{1994}

\bibitem{GalGW2019}
{\au C.G. Gal, M. Grasselli, H. Wu},
{\ti Global weak solutions to a diffuse interface model for incompressible two-phase flows with moving contact lines and different densities},
{\jou Arch.\ Ration.\ Mech.\ Anal.}
\no{234}{1--56}{2019}

\bibitem{GKL2018}
{\au M.H. Giga, A. Kirshtein, C. Liu},
{\bk Variational modeling and complex fluids},
{\jou Handbook of mathematical analysis in mechanics of viscous fluids},
\eds{73--113}{Springer, Cham}{ 2018}

\bibitem{G2020}
{\au A. Giorgini},
{\ti Well-posedness for a diffuse interface model for two-phase Hele-Shaw flows},
{\jou J.\ Math.\ Fluid Mech.}
\textbf{22}:5 (2020).

\bibitem{GGW2018}
{\au A. Giorgini, M. Grasselli, H. Wu},
{\ti The Cahn-Hilliard-Hele-Shaw system with singular potential},
{\jou Ann.\ Inst.\ H.\ Poincar\'{e} Anal.\ Non Lin\'{e}aire}
\no{35}{1079--1118}{2018}


\bibitem{GMT2019}
{\au A. Giorgini, A. Miranville, R. Temam},
{\ti Uniqueness and Regularity for the Navier-Stokes-Cahn-Hilliard system},
{\jou SIAM J.\ Math.\ Anal.}
\no{51}{2535--2574}{2019}

\bibitem{GT2020}
{\au A. Giorgini, R. Temam},
{\ti Weak and strong solutions to the nonhomogenenous incompressible Navier-Stokes-Cahn-Hilliard system},
{\jou J.\ Math.\ Pures Appl.}
(accepted for publication), 2020. 

\bibitem{GPV1996}
{\au M. E. Gurtin, D. Polignone, J. Vi\~{n}als},
{\ti Two-phase binary fluids and immiscible fluids described by an order parameter},
{\jou Math.\ Models Methods Appl. Sci.}
\no{6}{815--831}{1996}

\bibitem{HMR2012}
{\au M. Heida, J. M\'{a}lek, K.R. Rajagopal},
{\ti On the development and generalizations of Cahn–Hilliard equations within a thermodynamic framework},
{\jou Z.\ Angew.\ Math.\ Phys.}
\no{63}{145--169}{2012}

\bibitem{KZ2015}
{\au M. Kotschote, R. Zacher},
{\ti Strong solutions in the dynamical theory of compressible fluid mixtures},
{\jou Math.\ Models Meth.\ Appl.\ Sci.}
\no{25}{1217--1256}{2015}

\bibitem{LT1998}
{\au J. Lowengrub, L. Truskinovsky},
{\ti Quasi-incompressible Cahn-Hilliard fluids and topological transitions},
{\jou Proc.\ Roy.\ Soc.\ Lond.\ A}
\no{454}{2617--2654}{1998}

\bibitem{SSBZ2017}
{\au M. Shokrpour Roudbari, G. Şimşek, E.H. van Brummelen, K.G. van der Zee},
{\ti Diffuse-interface two-phase flow models with different densities: A new quasi-incompressible form and a linear energy-stable method},
{\jou Math.\ Models Meth.\ Appl.\ Sci.}
\no{28}{733--770}{2017}

\bibitem{Temam}
{\au R. Temam},
{\bk Navier-Stokes Equations and Nonlinear Functional Analysis}.
\eds{CBMS-NSF Regional Conference Series in Applied Mathematics, 66}{Society for Industrial and Applied Mathematics (SIAM), Philadelphia, PA}{1995}

\end{thebibliography}
\end{document}